\documentclass[opre,nonblindrev]{pinforms3} 

\OneAndAHalfSpacedXI 

\usepackage{algorithm}
\usepackage{algpseudocode}
\newfloat{algorithm}{t}{lop}
\newcommand{\Var}{\text{Var}}
\DeclareMathOperator{\diam}{diam}
\usepackage{xcolor}
\usepackage{setspace}

\usepackage{natbib}
\bibpunct[, ]{(}{)}{,}{a}{}{,}%
\def\bibfont{\small}%
\TheoremsNumberedThrough     
\ECRepeatTheorems

\EquationsNumberedThrough    

\begin{document}
\RUNAUTHOR{Mintz et al.}

\RUNTITLE{ROGUE Bandits}

\TITLE{Non-Stationary Bandits with Habituation and Recovery Dynamics}

\ARTICLEAUTHORS{%
\AUTHOR{Yonatan Mintz}
\AFF{School of Industrial and Systems Engineering, Georgia Institute of Technology, Atlanta, GA 30332, \EMAIL{ymintz@gatech.edu}}

\AUTHOR{Anil Aswani, Philip Kaminsky}
\AFF{Department of Industrial Engineering and Operations Research, University of California, Berkeley, CA 94720, \EMAIL{\{aaswani,kaminsky,\}@berkeley.edu}}

\AUTHOR{Elena Flowers}
\AFF{Department of Physiological Nursing, School of Nursing, University of California, San Francisco, CA 94143, \EMAIL{	elena.flowers@ucsf.edu}}

\AUTHOR{Yoshimi Fukuoka}
\AFF{Department of Physiological Nursing \& Institute for Health \& Aging, School of Nursing, University of California, San Francisco, CA 94143, \EMAIL{	yoshimi.fukuoka@ucsf.edu}}}

\ABSTRACT{%
Many settings involve sequential decision-making where a set of actions can be chosen at each time step, each action provides a stochastic reward, and the  distribution for the reward provided by each action is initially unknown. However, frequent selection of a specific action may reduce the expected reward for that action, while abstaining from choosing an action may cause its expected reward to increase. Such non-stationary phenomena are observed in many real world settings such as personalized healthcare-adherence improving interventions and targeted online advertising. Though finding an optimal policy for general models with non-stationarity is PSPACE-complete, we propose and analyze a new class of models called ROGUE (Reducing or Gaining Unknown Efficacy) bandits, which we show in this paper can capture these phenomena and are amenable to the design of policies with provable properties. We first present a consistent maximum likelihood approach to estimate the parameters of these models, and conduct a statistical analysis to construct finite sample concentration bounds. Using this analysis, we develop and analyze two different algorithms for optimizing ROGUE models: an upper confidence bound algorithm (ROGUE-UCB) and an $\epsilon$-greedy algorithm ($\epsilon$-ROGUE). Our theoretical analysis shows that under proper conditions the ROGUE-UCB and $\epsilon$-ROGUE algorithms can achieve logarithmic in time regret, unlike existing algorithms which result in linear regret. We conclude with a numerical experiment using real world data from a personalized healthcare-adherence improving intervention to increase physical activity. In this intervention, the goal is to optimize the selection of messages (e.g., confidence increasing vs. knowledge increasing) to send to each individual each day to increase adherence and physical activity. Our results show that ROGUE-UCB and $\epsilon$-ROGUE perform better in terms of aggregated regret and average reward when compared to state of the art algorithms, and in the context of this intervention the use of ROGUE-UCB increases daily step counts by roughly 1,000 steps a day (about a half-mile more of walking) as compared to other algorithms in a simulation experiment. 
}%



\maketitle

%

\section{Introduction}
Multi-armed bandits are commonly used to model sequential decision-making in settings where there is a set of actions that can be chosen at each time step, each action provides a stochastic reward, and the distribution for the reward provided by each action is initially unknown.  The problem of constructing a policy for sequentially choosing actions in multi-armed bandits requires balancing \emph{exploration} versus \emph{exploitation}, the tradeoff between selecting what is believed to be the action that provides the best reward and choosing other actions to better learn about their underlying distributions.  Bandit models have been applied in a variety of healthcare settings \citep{thompson1933,wang2011,bastani201,schell2016}.  For instance, \cite{bastani201} considered the problem of selecting drugs to give to a patient from a set (where each drug is an action) in order to treat a specific disease (the reward is the improvement in patient health in response to the drug); the bandit policy asymptotically identifies the optimal drug for that particular patient.   Other common applications involve online advertising \citep{agrawal2013,johari2015}, where selecting an ad to show is an action and the reward is the total number (from a large population) of viewers who click on the ad, as well as in various supply chain settings \citep{afeche2013,ban2014,caro2007}.

However, most bandit models assume that the distribution for the reward provided by each action is constant over time.  This is a reasonable assumption in a large number of applications, such as the ones described above.  However, many applications involve actions that are applied to a single individual, where the rewards depend upon behavioral responses of the individual to the applied actions.  In these behavioral settings, the response to a particular action is not generally stationary.  Frequent selection of a particular action will lead to habituation to that action by the individual, and the reward for that action will decrease each time it is selected.  For example, repeatedly showing the same ad to a single individual may cause the ad to become less effective in soliciting a response from that individual.  Furthermore, another complimentary phenomenon can also occur; refraining for a period of time from showing a particular ad to a single individual may cause the ad to become more effective when reintroduced.

Most techniques for designing policies for decision-making for multi-armed bandits apply to the setting where the rewards for each action are stationary.  However, designing a policy without considering the non-stationarity of a system (when the system is in fact non-stationary) often leads to poor results in terms of maximizing rewards \citep{besbes2014,hartland2006} because policies eventually converge to a stationary policy. The problem of designing policies for bandit models with non-stationarity has been studied in specific settings, but approaches in the literature are either computationally intractable, or the settings analyzed are not flexible enough to capture the habituation and recovery phenomenon described above.  The aim of this paper is to propose a flexible bandit model that is able to effectively model habituation and recovery, and to present an approach for designing an effective policy for this bandit model.

\subsection{Literature Review}

Data-driven decision-making can be categorized into batch formulations and online formulations.  Batch formulations \citep{aswani2016,mintz2017,ban2014,ban2015,bertsimas2014} use a large amount of data to estimate a predictive model and then use this model for optimization.  Adaptation to new data occurs by reestimating the predictive model, which is done periodically after a specified amount of additional data is collected.  

On the other hand, online formulations involve constructing a policy that is updated every time a new data point is collected.  Bandit models are a particularly important example of online formulations, and there has been much work on constructing policies for stationary bandits.  Approaches for designing policies for stationary bandits include those using upper confidence bounds \citep{auer2002finite,chang2005,bastani201}, Thompson sampling \citep{thompson1933,russo2014,russo2016,agrawal2013}, Bayesian optimization \citep{frazier2016,xie2013,xie2016}, knowledge gradients \citep{ryzhov2011,ryzhov2012}, robust optimization \citep{kim2015}, and adversarial optimization \citep{auer2002,agrawal2014,koolen2014,koolen2015}.


Restless bandits are a notable class of bandit models that capture non-stationarity, because choosing any single action causes the rewards of potentially all the actions to change.  Though dynamic programming  \citep{liu2010indexability,whittle1988restless}, approximation algorithms \citep{guha2010}, and mathematical programming \citep{bertsimas1994,bertsimas2000,caro2007} have been proposed as tools for constructing policies in this setting, the problem of computing an optimal policy for restless bandits is  PSPACE-complete \citep{papadimitriou1999}, meaning that designing policies that are approximately optimal is difficult. These models more broadly fall under the class of Partially Observable Markov Decision Processes (POMDPs), a class of problems where the decision maker is interested in solving a Markov Decision Process (MDP) with limited information about the system. The canonical methods for finding optimal policies for POMDPs involve the conversion of the POMDP into an equivalent belief MDP \citep{bertsekas2005} in which the states correspond to a probability distribution reflecting the decision maker's belief that they are in a specific state. While for small POMDPs it is possible to solve the belief MDP using dynamic programming approaches \citep{monahan1982,eagle1984,ayer2012}, in general, since POMDPs are PSPACE-complete, the resulting MDP is often too large to optimize effectively and requires using approximate dynamic programming \citep{ccelik2015,chen2018,deng2014,qi2017,shani2013,smith2012,wang2015}. Therefore, in recent literature, several solutions approaches have been proposed for solving POMDPs, which often involve exploiting the specific structure of the model in question. The framework of ROGUE multiarmed bandits we propose in this paper can be seen as a particular type of POMDP or restless bandit model, with a known deterministic transition model and partially observed rewards. To achieve efficient theoretical bounds and practical performance, the methods we propose are specific to this particular POMDP structure.



Another related research stream designs policies for non-stationary multi-armed bandits with specific structures.  For instance, model-free approaches have been proposed \citep{besbes2014,besbes2015non,garivier2008,anantharam1987} for settings with bounded variations, so that  rewards of each action are assumed to change abruptly but infrequently.  These policies have been shown to achieve $\mathcal{O}(\sqrt{T\log T})$ suboptimality in worst case scenarios but can achieve $\mathcal{O}(\log T)$ under proper model assumptions. Recently, there has been interest in studying more structured non-stationary bandits. Two relevant examples are Adjusted Upper Confidence Bounds (A-UCB) and rotting bandits \citep{bouneffouf2016,levine2017}, where each action has a set of unknown but stationary parameters and a set of known non-stationary parameters that characterize its reward distribution.  Policies designed for these settings achieve $\mathcal{O}(\log T)$ suboptimality, but these settings are unable to capture the habituation and recovery phenomenon that is of interest to us.


\subsection{ROGUE Bandits}

In this paper, we define the ROGUE (reducing or gaining unknown efficacy) bandit model, which can capture habituation and recovery phenomenon, and then we design a nearly-optimal policy for this model.  ROGUE bandits are appropriate for application domains where habituation and recovery are important factors for system design; we present two such examples, in online advertising and personalized healthcare, below.

\subsubsection{Personalized Healthcare-Adherence Improving Interventions}
One hundred fifty minutes of moderate-intensity aerobic physical activity each week has been shown to reduce the risk of cardiovascular disease, other metabolic disorders, and certain types of cancers \citep{physical2008,friedenreich2010,sattelmair2011,lewis2017}.  However, maintaining this level of moderate intensity activity is challenging for most adults. As such, proper motivation through providing daily exercise goals and encouragement has been found to be effective in helping patients succeed in being active \citep{fukuoka2011,fukuoka2014,fukuoka2014b}.

In recent years, there has been an increased rate of adoption of fitness applications and wearable activity trackers, making it easier and less costly to implement physical activity programs \citep{pwc2014}. These trackers and mobile applications record daily activity, communicate activity goals, and send motivational messages. Despite these digital devices having collected a large amount of personal physical activity data, many of the most popular activity trackers provide static and non-personalized activity goals and messages to their users \citep{rosenbaum2016}. Furthermore, the choice of motivational messages sent to users may have significant impact on physical activity, because if users receive similar messages too frequently they may become habituated and not respond with increased activity, while seldom sent messages may better increase activity due to their novelty and diversity.  Because the ROGUE bandits can model habituation and recovery of rewards for different actions, we believe they present a useful framework for the design of policies that choose which messages to send to users based on data consisting of what messages they received each day and the corresponding amounts of physical activity on those days.

Personalized healthcare has been extensively studied in the operations literature. \cite{aswani2016,mintz2017} explore the use of behavioral analytics to personalize diet and exercise goals for clinically supervised weight loss interventions in an offline setting. Markov decision processes have also been used for decision-making in personalized healthcare \citep{ayer2015, mason2013, deo2013, kucukyazici2011,leff1986,wang2011, gupta2008b,smilowitz2016,schell2016}. In contrast to bandit models where only the reward for the prescribed action can be observed, these methods broadly assume that the full state of the system can be observed, and thus do not require statistical estimation. Additionally, various multi-armed bandit approaches \citep{bastani201,wang2011} have also been proposed for healthcare problems where habituation and recovery are not significant factors. 

\subsubsection{Online Content Creation and Advertising}
Online advertising is one of the fastest-growing industries in the US. In fact, as of 2016, US Internet advertising spending has increased to over \$72.5 billion, surpassing the amount spent on TV ads \citep{richter2017}. However, as this form of advertising becomes more prevalent, advertisers have been struggling to ensure that ads retain there effectiveness.
This has been attributed to Internet users being habituated by impersonal and standardized ads \citep{goldfarb2014,portnoy2010} which are rarely varied. For these reasons, there has been significant interest in the operations literature in creating automated systems that can utilize user-level data to better target and customize ads \citep{ghose2009,goldfarb2011}. In particular, since the effect of a no-longer-effective advertisement may recover after a user has not seen it for some period of time, incorporating recovery and habituation dynamics into advertising models could yield more effective advertising campaigns.

In general, multi-armed bandit models have been proposed to model online advertising, where each action corresponds to a different type of advertisement, and the reward is equivalent to either a conversion or a click from a prospective consumer. Several approaches have been used to design such an ad targeting system, including adversarial and stochastic multi-armed bandit models \citep{bertsimas2007,chen2013,kleinberg2008,liu2010,yi2010}, and online statistical testing \citep{johari2015}. However, while some of these approaches use contextual data to better serve ads to individuals, they are still designed under assumptions of stationarity.  As a result, these approaches will lead to policies that show duplicated ads to individuals, which can potentially causing habituation, whereas other ads that might have recovered efficacy may not be served at all. In contrast, ROGUE Bandit models can explicitly consider the time-varying efficacy each type of ad, and thus directly capture user habituation to a specific ad, and track the recovery of efficacy of a particular ad for a specific individual.

\subsection{Outline}
In Section \ref{sec:mab}, we formally introduce the ROGUE bandit model. To the best of our knowledge, this is the first work where a non-stationary bandit model has been defined that is able to capture habituation and recovery phenomenon, and is at the same time amenable to the design of nearly-optimal policies.  Because the ROGUE bandit is a general model, we describe two specific instantiations: the ROGUE generalized linear model and the ROGUE agent. 
 
Next, in Section \ref{sec:parm_est} we analyze the problem of estimating the parameters of a single action.  We present a statistical analysis of maximum likelihood estimation (MLE) for a single action, and use empirical process theory to derive finite sample bounds for the convergence of parameters estimates.  Specifically, we show that the MLE estimates converge to the true parameters at a $1/\sqrt{T}$ rate.

Section \ref{sec:ducb} describes an upper-confidence bound policy for ROGUE bandits, and we call this policy the ROGUE-UCB algorithm.  The main result of this section is a rigorous $\mathcal{O}(\log T)$ bound on the suboptimality of the policy in terms of regret, the difference between the reward achieved by the policy and the reward achieved by an optimal policy.  Our $\mathcal{O}(\log T)$ bound is significant because this is the optimal rate achievable for approximate policies in the stationary case \citep{lai1985}.  We prove our bound using methods from the theory of concentration of measure.

We conclude with Section \ref{sec:num_exp}, where we introduce a ``tuned'' version of ROGUE-UCB and then conduct numerical experiments to compare the efficacy of our ROGUE-UCB algorithm to other policies that have been developed for bandit models.  Our experiments involve two instantiations of ROGUE bandit models.  First, we compare different bandit policies using a ROGUE generalized linear bandit to generate data.  Second, we compare different bandit policies using a ROGUE agent to generate data, where the parameters of this bandit model are generated using data from a physical activity and weight loss clinical trial \citep{fukuoka2014}.  This second experiment specifically addresses the question of how to choose an optimal sequence of messages to send to a particular user in order to optimally encourage the user to increase physical activity, and it can be interpreted as a healthcare-adherence improving intervention.  Our experiments show that ROGUE-UCB outperforms all other considered bandit policies, and that it achieves logarithmic regret, in contrast to other bandit algorithms that achieve linear regret.

\section{Defining Reducing or Gaining Unknown Efficacy (ROGUE) Bandits}
\label{sec:mab}

This section first describes the stationary multi-armed bandit (MAB) model, in order to emphasize modeling differences in comparison to our ROGUE bandit model that is introduced in this section.  Our goal in defining ROGUE bandits is to have a model that can capture specific non-stationary phenomena found in behavioral applications, and so we next formally introduce the model elements of ROGUE bandits.  To provide better intuition about ROGUE bandits, we also present two specific instantiations of a ROGUE bandit that incorporate different behavioral effects.

\subsection{Stationary MAB Model}

The stationary MAB is a setting where there is a finite set of actions $\mathcal{A}$ that can be chosen at each time step $t$, each action $a \in \mathcal{A}$ provides a stochastic reward $r_a$ with distribution $\mathbb{P}_{\theta_a}$, and the parameters $\theta_a\in\Theta$ for $a \in\mathcal{A}$ are constants that are initially unknown but lie in a known compact set $\Theta$.  The problem is to construct a policy for sequentially choosing actions in order to maximize the expected reward.  More specifically, let $\pi_t\in\mathcal{A}$ be the action chosen at time $t = 1,\ldots,T$.  Then the policy consists of functions $\pi_t(r_{\pi_1},\ldots,r_{\pi_{t-1}},\pi_1,\ldots,\pi_{t-1})\in\mathcal{A}$ that depend on past rewards and actions.  For notational convenience, we will use $\Pi = \{\pi_t(\cdot)\}_{t=1}^T$ to refer to the policy.  In this notation, the problem of constructing an optimal policy to maximize expected reward can be written as $\max_{\Pi \in \mathcal{A}^T}\sum_{t=1}^T \mathbb{E}r_{\pi_t}$.  Note that certain regularity is needed from the distributions to ensure this maximization problem is well-posed.  One common set of assumptions is that the distributions $\mathbb{P}_{\theta_a}$ for $a\in\mathcal{A}$ are sub-Gaussian, and that the reward distributions are all independent.  

For the stationary MAB, we can define an optimal action $a^* \in \mathcal{A}$, which is any action such that $\mathbb{E}r_{a^*} \geq \mathbb{E}r_a$ for all $a\in\mathcal{A}$.  The benefit of this definition is it allows us to reframe the policy design problem in terms of minimizing the cumulative expected regret $\mathbb{E}R_\Pi(T) = \mathbb{E}[  Tr_{a^*} - \sum_{i=1}^{T}r_{\pi_t}]$, where the quantity $r_{a^*} - r_{\pi_t}$ is known as the regret at time $t$.  Observe that minimizing $\mathbb{E}R_\Pi(T)$ is equivalent to maximizing $\sum_{t=1}^T \mathbb{E}r_{\pi_t}$.  It has been shown by \cite{gittins1979} that an index policy is optimal for the stationary MAB. Since these indexing policies are difficult to compute, other approximate policies have been proposed \citep{lai1985,auer2002finite}. Some of the most common policies use upper confidence bounds \citep{auer2002finite,garivier2011}, which take actions optimistically based on estimates of the parameters $\theta_a$. Unfortunately, it has been shown that these index policies and upper confidence bound policies can have arbitrarily bad performance in a non-stationary setting \citep{hartland2006,besbes2014}.



\subsection{Reducing or Gaining Unknown Efficacy (ROGUE) Bandits} \label{sec:dmab}

A disadvantage of the stationary MAB is that it does not allow rewards to change over time in response to previous actions, and this prevents the stationary MAB model from being able to capture habituation or recovery phenomena.  Here, we define ROGUE bandits that can describe such behavior.  The ROGUE bandit is a setting where there is a finite set of actions $\mathcal{A}$ that can be chosen at each time step $t$, each action $a \in \mathcal{A}$ at time $t$ provides a stochastic reward $r_{a,t}$ that has a sub-Gaussian distribution $\mathbb{P}_{\theta_a,x_{a,t}}$ with expectation $\mathbb{E}r_{a,t} = g(\theta_a,x_{a,t})$ for a bounded function $g$, the parameters $\theta_a\in\Theta$ for $a \in\mathcal{A}$ are constants that are initially unknown but lie in a known compact, convex 
set $\Theta$, and each action $a\in\mathcal{A}$ has a state $x_{a,t}$ with nonlinear dynamics $x_{a,t+1} = h_a(x_{a,t},\pi_{a,t})$
where $\pi_{a,t} = \mathbf{1}[\pi_t = a]$,  $h:\mathcal{X}\times \mathbb{B} \ \rightarrow \mathcal{X}$ is a known dynamics function, and $\mathcal{X}$ is a compact convex set such that $x_{a,t} \in \mathcal{X} \ \forall a,t$, and $x_{a,0}$ is initially unknown for $a\in\mathcal{A}$.  This particular model formulation is useful for various real world scenarios in which effective system models exist but states are not directly observed due to system or measurement noise. This can be achieved by training models on previously collected controlled data, or by using models based on information from qualified domain experts. In the case of personalized healthcare-adherence interventions and on-line marketing, often times effective models can be generated from previous experimental data (e.g., from RCTs or focus group testing) though state observation in non-laboratory (i.e., not completely controlled) situations is challenging. This data can be used in conjunction with MLE techniques to obtain useful estimates of the system dynamics. In cases where limited data and expertise are available, alternative models ought to be used, the discussion of which is out of the scope of this paper.

The problem is to construct a policy for sequentially choosing actions in order to maximize the expected reward.  Observe that the ROGUE bandit model is non-stationary since the reward distributions depend upon previous actions. This makes the problem of designing policies more difficult than that of designing policies for the stationary MAB.  More specifically, let $\pi_t\in\mathcal{A}$ be the action chosen at time $t = 1,\ldots,T$.  Then the policy consists of functions $\pi_t(r_{\pi_1},\ldots,r_{\pi_{t-1}},\pi_1,\ldots,\pi_{t-1})\in\mathcal{A}$ that depend on past rewards and actions.  For notational convenience, we will use $\Pi = \{\pi_t(\cdot)\}_{t=1}^T$ to refer to the policy.  In this notation, the problem of constructing an optimal policy to maximize expected reward can be written as
\begin{equation}
\textstyle\max_{\Pi \in \mathcal{A}^T} \{ \sum_{t=1}^T g(\theta_{\pi_t},x_{\pi_t,t}) : x_{a,t+1} = h_a(x_{a,t},\pi_{a,t}) \text{ for } a \in \mathcal{A},\ t \in \{0,...,T-1\}\}.
\end{equation}
This can be reframed as minimizing expected cumulative regret \citep{besbes2014,garivier2008,bouneffouf2016}: Unlike the stationary MAB, we cannot define an optimal action, but rather must define an optimal policy $\Pi^* = \{\pi^*_t(\cdot)\}_{t=0}^{T}$, which can be thought of as an oracle that chooses the optimal action at each time step.  Then the problem of designing an optimal policy is equivalent to minimizing $R_\Pi(T) = \sum_{t=1}^T r_{\pi^*_t,t} - r_{\pi_t,t}$ or in expectation $\mathbb{E}R_\Pi(T) = \sum_{t=1}^T g(\theta_{\pi^*_t},x_{\pi^*_t,t}) - g(\theta_{\pi_t},x_{\pi_t,t})$ subject to the state dynamics defined above, a notion known as either dynamic or tracking regret. Note that since the optimal policy may be quite different from the policy actually implemented, the resulting states and rewards observed at each time period may not necessarily be the same. There have been stronger definitions of dynamic regret proposed in the literature including strongly adaptive regret \citep{daniely2015,zhang2017}. These definitions have been shown to not hold generally in the bandit feedback case even for models with optimal tracking regret rates \citep{daniely2015}, and are applied in the context of full information convex on-line optimization with limited total variation. In this paper, we center our analysis on tracking regret as our setting includes both bandit feedback and is model-based without necessarily bounded total variation.

\subsection{Technical Assumptions on ROGUE Bandits}

In this paper, we will design a policy for ROGUE bandits that follow the assumptions described below:

\begin{assumption} \label{ass:indep_ass}
The rewards $r_{a,t}$ are conditionally independent given $x_{a,0},\theta_a$ (or equivalently the complete sequence of $x_{a,t}, \pi_t$ and $\theta_a$). 
\end{assumption}
This assumption states that for any two time points $t,t'$ such that $t \neq t'$ we have that $r_{a,t}|\{x_{a,t},\theta\} $ is independent of $r_{a,t'}|\{x_{a,t'},\theta\}$, and it is a mild assumption because it is the closest analogue to the assumption of independence of rewards in the stationary MAB.

\begin{assumption} \label{ass:log_conv_ass}
The reward distribution $\mathbb{P}_{\theta,x}$ has a log-concave probability density function (p.d.f.) $p(r|\theta,x)$ for all $x\in\mathcal{X}$ and $\theta\in\Theta$.
\end{assumption}
This assumption provides regularity for the reward distributions, and is met by many common distributions (e.g., Gaussian and Bernoulli).

Now define $f(\cdot)$ to be $L$-Lipschitz continuous if $|f(x_1) - f(x_2)| \leq L\|x_1-x_2\|_2$ for all $x_1,x_2$ in the domain of $f$. We define $f$ to be locally $L$-Lipschitz on a compact set $\mathcal{S}$ if it has the Lipschitz property for all points of its domain on set $\mathcal{S}$. Our next assumption is on the stability of the above distributions with respect to various parameters.
\begin{assumption}\label{ass:lip_ass}
The log-likelihood ratio $\ell(r;\theta',x',\theta,x) = \log\frac{p(r|\theta',x')}{p(r|\theta,x)}$ associated with the distribution family $\mathbb{P}_{\theta,x}$ is locally $L_f$-Lipschitz continuous with respect to $x,\theta$  on the compact set $\mathcal{X}\times\Theta$ for all values of $\theta',x'\in \mathcal{X}\times\Theta$, and $g$ is locally $L_g$-Lipschitz continuous with respect to $x,\theta$ on the compact set $\mathcal{X}\times\Theta$.
\end{assumption} 

This assumption ensures that if two sets of parameters are close to each other in value then the resulting distributions will also be similar. We make the following additional assumption about the functional structure of the reward distribution family:

\begin{assumption} \label{ass:sub_gauss}
The reward distribution $\mathbb{P}_{\theta,x}$ for all $\theta\in\Theta$ and $x\in\mathcal{X}$ is sub-Gaussian with parameter $\sigma$, and either $p(r|\theta,x)$ has a finite support or $\ell(r;\theta',x',\theta,x)$ is locally $L_p$-Lipschitz with respect to $r$. 
\end{assumption} 

This assumption (or a similar type of regularity) is needed to ensure that sample averages are close to their means, and it is satisfied by many distributions (e.g., a Gaussian location family with known variance). 

Last, we impose conditions on the dynamics for the state of each action:

\begin{assumption} \label{ass:h_lip}
The dynamic transition function $h$ is $L_h$ Lipschitz continuous such that $L_h \leq 1$.
\end{assumption}

This assumption is needed to ensure the states of each action do not change too quickly, and it implies that the dynamics are stable. 

\subsection{Instantiations of ROGUE Bandits}

The above assumptions are general and apply to many instantiations of ROGUE bandit models. To demonstrate the generality of these assumptions, we present two particular instances of ROGUE bandit models.

\subsubsection{ROGUE Agent} \label{sec:stack_band}

Our first instantiation of a ROGUE bandit model consists of a dynamic version of a principal-agent model \citep{stackelberg1952,radner1985,laffont2002,mintz2017}, which is a model where a principal designs incentives to offer to an agent who is maximizing an (initally unknown to the principal) utility function that depends on the incentives.  In particular, consider a setting with a single (myopic) agent to whom we would like to assign a sequence of behavioral incentives $\pi_t \in \mathcal{A}$, and the states $x_{a,t}$ and parameters $\theta_a$ are scalars.  Given a particular incentive $\pi_t$ at time $t$, the agent responds by maximizing the (random) utility function
\begin{equation}
\label{eqn:ragent}
\textstyle r_t = \argmax_{r \in [0,1]} -\frac{1}{2}r^2 - (c_{a,t}+ \sum_{a \in \mathcal{A}}x_{a,t}\pi_{t,a})r,
\end{equation}
where for fixed $a\in\mathcal{A}$ we have that $c_{a,t}$ are i.i.d. random variables with a distribution $\mathbb{P}_{\theta_a}$ such that $\mathrm{Var}(c_{a,t}) = \sigma^2(\theta_a) < \infty$ and $\sigma^2: \mathbb{R}\rightarrow \mathbb{R}_+$ is invertible. Moreover, the state dynamics are
\begin{equation}
x_{a,t+1} = \text{proj}_\mathcal{X}(\alpha_a x_{a,t} + b_a (1- \pi_{a,t}) -k_a),
\end{equation}
  Note the distribution of $r_{t}$ is fully determined by $x_{a,t},\theta_a,\{\pi_k\}_{k=0}^t$, which means the rewards satisfy Assumption \ref{ass:indep_ass}.

We can further analyze the above ROGUE agent model.  Solving the agent's optimization problem (\ref{eqn:ragent}) gives
 \begin{equation}
 r_t|\{x_{a,t},\theta_a\} = \begin{cases} 0 & \text{if } \; c_{a,t} \leq - x_{a,t}, \\
 1 & \text{if } \; c_{a,t} \geq 1- x_{a,t}, \\
  c_{a,t} + x_{a,t} & \text{otherwise}
 \end{cases}
 \end{equation}
 We can express the distribution of $r_t|\{x_{a,t},\theta_a\}$ in terms of the cumulative distribution function (c.d.f.) $F(\cdot)$ and p.d.f. $f(\cdot)$ of $c_{a,t}$:
\begin{equation}
p(r_t|\{x_{a,t},\theta_a\}) = F(-x_{a,t})\delta(r_t) +(1-F(1-x_{a,t})) \delta(1-r_t) + f(r_t -x_{a,t}) \mathbf{1}[r_t \in (0,1)].
\end{equation}
Though $p(r_t|\{x_{a,t},\theta_a\})$ is not an absolutely continuous function, it satisfies Assumptions \ref{ass:log_conv_ass} and \ref{ass:lip_ass}, whenever $c_t$ has a log-concave p.d.f. that is Lipschitz continuous, if we interpret the above probability measure $p(r_t|\{x_{a,t},\theta_a\})$ as a p.d.f.


\subsubsection{ROGUE Generalized Linear Model (GLM)}
\label{sec:dglm}
Dynamic logistic models and other dynamic generalized linear models \citep{mccullagh1984,filippi2010} can be interpreted as non-stationary generalizations of the classical (Bernoulli reward) stationary MAB \citep{gittins1979,lai1985,garivier2011}. Here, we further generalize these models: Consider a setting where $r_{a,t}|\{\theta_a,x_{a,t}\}$ is an exponential family with mean parameter
\begin{equation}
\mu_{a,t} = \mathbb{E}r_t = g(\alpha_a^T\theta_a + \beta_a^Tx_{a,t}),
\end{equation}
for known vectors $\alpha_a,\beta_a$, where the action states $x_{a,t}$ have appropriate dynamics.  In this situation, we can interpret $g(\cdot)$ as a link function of a generalized linear model (GLM). For example, if $g$ is a logit function, then this model implies the rewards have a Bernoulli distribution with parameter
\begin{equation}
\mu_{a,t} = \frac{1}{1+\exp(-(\alpha_a^T\theta_a + \beta_a^Tx_{a,t}))}.
\end{equation}
For the logistic case, the $r_{a,t}$ is bounded and satisfies Assumptions \ref{ass:indep_ass}-\ref{ass:log_conv_ass}. These assumptions are also satisfied if $r_{a,t}$ can be linked to a truncated exponential family distribution restricted to $[0,1]$, meaning if the p.d.f. of $r_{a,t}|\{x_{a,t},\theta_a \}$ is 
\begin{equation}
\frac{h(r)}{F(1)-F(0)} \exp\big(T(r)g(\alpha_a^T\theta_a + \beta_a^Tx_{a,t}) - A(\alpha_a^T\theta_a + \beta_a^Tx_{a,t})\big),
\end{equation}
where $T(r)$ is a sufficient statistic $h$ is the base measure (not to be confused with dynamics $h_a$), $A$ is the log-partition function, and $F$ is the full cdf of the un-truncated distribution.  If instead we consider sub-Gaussian exponential families with infinite support, Assumption \ref{ass:sub_gauss} is satisfied if the sufficient statistic of the GLM is Lipschitz or bounded with respect to $r$. While we will mainly consider one-dimensional rewards (i.e., $r_{a,t}\in\mathbb{R}$), we note that this framework can also be extended to vector and array dynamic GLM's.

\section{Parameter Estimation for ROGUE Bandits}
\label{sec:parm_est}
Our approach to designing a policy for ROGUE bandits will involve generalizing the upper confidence bound policies \citep{auer2002finite,chang2005,bastani201} that have been developed for variants of stationary MAB's.  As per the name of these policies, the key step involves constructing a confidence bound for the parameters $\theta_a, x_{a,0}$ characterizing the distribution of each action $a\in\mathcal{A}$.  This construction is simpler in the stationary case because the i.i.d. structure of the rewards allows use of standard Chernoff-Hoeffding bounds \citep{wainwright2015}, but we can no longer rely upon such i.i.d. structure for ROGUE bandits which are fundamentally non-stationary. This is because in ROGUE bandits the reward distributions depend upon states $x_{a,t}$, and so the structure of ROGUE bandits necessitates new theoretical results on concentration of measure in order to construct upper confidence bounds for the relevant parameters.

For this analysis, let the variables $\{r_{a,t}\}_{t=1}^T$ be the observed rewards for action $a\in\mathcal{A}$. It is important to note that the $r_{a,t}$ here are no longer random variables, but are rather the actual observed values.  Since the reward distributions for each action are mutually independent by the dynamics $h_a$, we can study the estimation problem for only a single action. Specifically, consider the likelihood $p( \{r_{a,t}\}_{t\in \mathcal{T}_a}| \theta_{a},x_{a,0} )$, where $\mathcal{T}_a \subset \{1,...,T\}$ is the set of times when action $a$ was chosen (i.e., $\pi_t = a$ for $t\in\mathcal{T}_a$).  Let $n(\mathcal{T}_a)$ denote the cardinality of the set $\mathcal{T}_a$.  Using Assumption \ref{ass:indep_ass}, the likelihood can be expressed as
\begin{equation}
\label{eqn:objmle}
p(\{r_{a,t}\}_{t\in \mathcal{T}_a}|\theta_a,x_{a,0}) = \prod_{t\in \mathcal{T}_a}p(r_{a,t}|\theta_a,x_{a,t}) \prod_{t\in \mathcal{T}_a}p(x_t|\theta_a,x_{a,t_-}).
\end{equation}
where $t_- = \max\{s\in\mathcal{T}_a : s < t\}$ is the latest observation before time $t$. Note the MLE of $\theta_a,x_{a,0}$ is $(\hat{\theta}_a,\hat{x}_{a,0}) \in \argmax \prod_{t\in \mathcal{T}_a}p(r_{a,t}|\theta_a,x_{a,t}) \prod_{t\in \mathcal{T}_a}p(x_t|\theta_a,x_{a,t_-})$. Observe that since the dynamics $h_a$ are known the one step likelihood $p(x_t|\theta_a,x_{a,t-1})$ is a degenerate distribution with all probability mass at $x_{a,t}$, by perpetuation of the dynamics $h$ with initial conditions $x_{a,t-1}$. Thus we can express the MLE as the solution to the constrained optimization problem
\begin{equation}
\textstyle (\hat{\theta}_a,\hat{x}_{a,0}) = \arg\min \{-\sum_{t\in \mathcal{T}_a}\log p(r_{a,t}|\theta_a,x_{a,t}) : x_{a,t+1} = h_a(x_{a,t},\pi_{a,t}) \text{ for } t \in \{0,\ldots,T\} \},
\end{equation}
where we have also taken the negative logarithm of the likelihood (\ref{eqn:objmle}). In this section, we will consider concentration properties of the solution to the above optimization problem. In particular, since our results concern the joint distributions of a trajectory of values, we require the definition of the following quantity which we call the trajectory Kullback--Leibler (KL) Divergence:

	\begin{definition}
	For some input action sequence $\pi_1^T$ and arm $a \in \mathcal{A}$ with dynamics $h_a$, given starting parameter values $(\theta_a,x_{a,0}), (\theta'_a,x'_{a,0}) \in \mathcal{X}\times\Theta$, define the trajectory KL--Divergence between these two trajectories with the same input sequence and different starting conditions as:
	\begin{equation}
	D_{a,\pi_1^T}(\theta_a,x_{a,0}||\theta'_a,x'_{a,0}) = \sum_{t\in \mathcal{T}_a}D_{KL}(\mathbb{P}_{\theta_a,x_{a,t}}||  \mathbb{P}_{\theta'_a,x'_{a,t}}) =\sum_{t\in \mathcal{T}_a}D_{KL}(\mathbb{P}_{\theta_a,h_a^t(x_{a,0})}|| \mathbb{P}_{\theta'_a,h_a^t(x'_{a,0})})
	\end{equation}
	where $h^k_a$ represents the functional composition of $h_a$ with itself $k$ times subject to the given input sequence, $\mathbb{P}_{\theta,x}$ is the probability law of the system under parameters $\theta,x$, and $D_{KL}$ is the standard KL-Divergence.
	\end{definition}

 If $\theta^*_a,x^*_{0,a}$ for $a\in\mathcal{A}$ are the true parameter values of a ROGUE Bandit model, then we show that 
\begin{theorem} \label{thm:concent_eq}
For any constant $\xi >0$ we have
\begin{equation}
P\Bigg(\frac{1}{n(\mathcal{T}_a)}D_{a,\pi_1^{T}}(\theta^*_a,x^*_{a,0}||\hat{\theta}_a,\hat{x}_{a,0}) \leq \xi +  \frac{c_f(d_x,d_\theta)}{\sqrt{n(\mathcal{T}_a)}}\Bigg) \geq 1 - \exp\Bigg(\frac{-\xi^2 n(\mathcal{T}_a)}{2 L_p^2 \sigma^2}\Bigg)
\end{equation}
where
\begin{equation}
c_f(d_x,d_\theta) = 8L_f\diam(\mathcal{X})\sqrt{\pi} + 48\sqrt{2}(2)^\frac{1}{d_x+d_\theta}L_f\diam(\mathcal{X}\times\Theta)\sqrt{\pi(d_x+d_\theta)}
\end{equation}
is a constant that depends upon $d_x$ (the dimensionality of $\mathcal{X}$) and $d_\theta$ (the dimensionality of $\Theta$), and $D_{a,\pi_1^T}(\theta_a,x_{a,0}||\theta'_a,x'_{a,0})$ is the trajectory Kullback--Leibler (KL) divergence between two different initial conditions.
\end{theorem}

\subsection{Conceptual Reformulation of MLE}
Our analysis begins with a reformulation of the MLE that removes the constraints corresponding to the dynamics through repeated composition of the function $h_a$.
\begin{proposition}\label{prop:reform_prop}
Let $\theta^*_a\in\Theta$ and $x^*_{a,0} \in \mathcal{X}$ for $a\in\mathcal{A}$ be the true underlying parameters of the system, then the MLE is given by 
\begin{equation}
\begin{aligned}
(\hat{\theta}_a,\hat{x}_{a,0}) = \argmin_{\theta_a,x_{a,0}\in \Theta \times \mathcal{X}} \frac{1}{n(\mathcal{T}_a)}\sum_{t\in \mathcal{T}_a}\log\frac{ p(r_{a,t}|\theta^*_a,h_a^t(x^*_{a,0},\theta^*_{a},\pi_1^t))}{p(r_{a,t}|\theta_a,h_a^t(x_{a,0},\theta_{a},\pi_1^t))}
\end{aligned}
\end{equation}
where the notation $h_a^k$ represents the repeated functional composition of $h_a$ with itself $k$ times, and $\pi_1^t$ is the sequence of input decisions from time 1 to time $t$.
\end{proposition}
The complete proof for this proposition is found in Appendix \ref{ap:proof}, and here we provide a sketch of the proof. Observe that this formulation is obtained by first adding constant terms equal to the likelihood of the true parameter values to the objective function and dividing by the total number of observations (which does not change the optimal solution), and then composing our system dynamics and writing them as explicit functions of the initial conditions. In practice, this reformulation is not practical to solve since clearly $\theta^*_a,x^*_{a,0}$ are not known \emph{a priori} and the composite function $h_a^t$ may have a complex form. However, for theoretical analysis this reformulation is quite useful, since for fixed $\theta_a,x_{a,0}$ taking the expected value of the objective under $\mathbb{P}_{\theta^*_a,x^*_{a,0}}$ yields
\begin{multline}
\mathbb{E}_{\theta^*_a,x^*_{a,0}}\frac{1}{n(\mathcal{T}_a)} \sum_{t\in\mathcal{T}_a}\log\frac{ p(r_{a,t}|\theta^*_a,h_a^t(x^*_{a,0},\theta^*_{a},\pi_1^t))}{p(r_{a,t}|\theta_a,h_a^t(x_{a,0},\theta_{a},\pi_1^t))} = \frac{1}{n(\mathcal{T}_a)} \sum_{t\in \mathcal{T}_a}D_{KL}(\mathbb{P}_{\theta^*_a,x^*_{a,t}}|| \mathbb{P}_{\theta_a,x_{a,t}})\\ =\frac{1}{n(\mathcal{T}_a)}D_{a,\pi_1^T}(\theta^*_a,x^*_{a,0}||\theta_a,x_{a,0}).
\end{multline} 
Essentially, we have reformulated the MLE problem in terms of minimizing the KL divergence between the trajectory distribution of potential sets of parameters to the trajectory distribution of the true parameter set. Since we have clear interpretation for the expectation of our objective function we can now proceed to compute concentration inequalities.
\subsection{Uniform Law of Large Numbers for ROGUE Bandits}

Since our estimates are computed by solving an optimization problem, a pointwise law of large numbers is insufficient for our purposes since such a result would not be strong enough to imply convergence of the optimal solutions.  To obtain proper concentration inequalities we must consider a uniform law of large numbers for the MLE problem.
\begin{theorem} \label{thm:two_side_ineq}
For any constant $\xi > 0$ we have
\begin{multline}
P\Bigg(\sup_{\theta_a,x_{a,0} \in \Theta\times \mathcal{X}}\Bigg|\frac{1}{n(\mathcal{T}_a)}\sum_{t\in\mathcal{T}_a}\log\frac{ p(r_{a,t}|\theta^*_a,h_a^t(x^*_{a,0},\theta^*_{a},\pi_1^t))}{p(r_{a,t}|\theta_a,h_a^t(x_{a,0},\theta_{a},\pi_1^t))}\\ - \frac{1}{n(\mathcal{T}_a)}D_{a,\pi_1^T}(\theta^*_a,x^*_{a,0}||\theta_a,x_{a,0}) \Bigg| > \xi +   \frac{c_f(d_x,d_\theta)}{\sqrt{n(\mathcal{T}_a)}}\Bigg) \leq \exp\Bigg(\frac{-\xi^2 n(\mathcal{T}_a)}{2L_p^2\sigma^2}\Bigg)
\end{multline}
where
\begin{equation}
c_f(d_x,d_\theta) = 8L_f\diam(\mathcal{X})\sqrt{\pi} + 48\sqrt{2}(2)^\frac{1}{d_x+d_\theta}L_f\diam(\mathcal{X}\times\Theta)\sqrt{\pi(d_x+d_\theta)}
\end{equation}
is a constant.
\end{theorem}

We will prove this result in several steps, the first of which uses the following lemma:

\begin{lemma} \label{lemma:lip_lem}
Consider the mapping
\begin{equation}
\varphi\Big(\{r_t\}_{t=1}^{n(\mathcal{T}_a)}\Big)=\sup_{\theta_a,x_{a,0} \in \Theta\times \mathcal{X}}\Bigg|\frac{1}{n(\mathcal{T}_a)}\sum_{t\in\mathcal{T}_a}\log\frac{ p(r_{a,t}|\theta^*_a,h_a^t(x^*_{a,0},\theta^*_{a},\pi_1^t))}{p(r_{a,t}|\theta_a,h_a^t(x_{a,0},\theta_{a},\pi_1^t))} - \frac{1}{n(\mathcal{T}_a)}D_{a,\pi_1^T}(\theta^*_a,x^*_{a,0}||\theta_a,x_{a,0}) \Bigg|.
\end{equation}
The mapping $\varphi$ is $L_p$-Lipschitz with respect to $\{r_t\}_{t=1}^{n(\mathcal{T}_a)}$.
\end{lemma}
A detailed proof is provided in Appendix \ref{ap:proof}, and the main argument of the proof relies on the preservation of Lipschitz continuity through functional composition and pointwise maximization. This result is necessary since showing that objective value variations are bounded is a prerequisite for the formalization of concentration bounds. Next we consider the Lipschitz constant of the log-likelihood with respect to the parameters.

\begin{lemma} \label{lem:lip_param}
For any $r\in\mathbb{R}$, $\bar{\theta}\in\Theta$, $\bar{x} \in\mathcal{X}$, define the function $\ell: \Theta\times\mathcal{X}\times\{1,...,T\} \rightarrow \mathbb{R}$ such that $\ell(\theta,x,t) = \log\frac{p(r|\bar{\theta},h_a^t(\bar{x},\bar{\theta},\pi_1^t))}{p(r|\theta,h_a^t(x,\theta,\pi_1^t))}$.  Then for fixed $t$, the function $\ell$ is Lipshitz with constant $L_f$. Moreover, for all $(x,\theta) \in \mathcal{X}\times\Theta$ and for all $t,t' \in \{1,...,T\}$ we have that $|\ell(\theta,x,t) -\ell(\theta,x,t')| \leq L_f\diam(\mathcal{X})$, where $\diam(\mathcal{X}) = \max_{x\in \mathcal{X}} \|x\|_2$.
\end{lemma}
The result of this lemma can be derived using a similar argument to that of Lemma \ref{lemma:lip_lem}, by noting that the dynamics are bounded and Lipschitz, and then applying Assumption \ref{ass:lip_ass}. The full proof of this lemma is in Appendix \ref{ap:proof}. Next we show the expected behavior of $\pi$ is bounded.

\begin{lemma} \label{lem:rademach}
Let $\varphi$ be defined as in Lemma \ref{lemma:lip_lem}. Then $\mathbb{E}\varphi(\{r_t\}_{t=1}^{n(\mathcal{T}_a)}) \leq \frac{c_f(d_x,d_\theta)}{\sqrt{n(\mathcal{T}_a)}}$, where 
\begin{equation}
c_f(d_x,d_\theta) = 8L_f\diam(\mathcal{X})\sqrt{\pi} + 48\sqrt{2}(2)^\frac{1}{d_x+d_\theta}L_f\diam(\mathcal{X}\times\Theta)\sqrt{\pi(d_x+d_\theta)}.
\end{equation}
\end{lemma}
The result of this lemma is derived by first using a symmetrization argument to bound the expectation by a Rademacher average and then using metric entropy bounds to derive the final result, and a complete proof is found in Appendix \ref{ap:proof}. Additional insight into these results is provided by the following remarks:
\begin{remark}
The result of Lemma \ref{lem:rademach} implies that $\mathbb{E}\varphi(\{r_{a,t}\}_{t=1}^{n(\mathcal{T}_a)}) = \mathcal{O}(\sqrt{\frac{d_x+d_\theta}{n(\mathcal{T}_a)}})$
\end{remark}
\begin{remark}
An improved constant can be achieved by using weaker metric entropy bounds (namely the union bound) however this would yield a bound of order $\mathcal{O}(\sqrt{\frac{(d_x+d_\theta)\log n(\mathcal{T}_a)}{n(\mathcal{T}_a)}})$
\end{remark}
Using the results of Lemmas \ref{lemma:lip_lem}--\ref{lem:rademach}, we can complete the sketch of the proof for Theorem \ref{thm:two_side_ineq}. Lemma \ref{lemma:lip_lem} says the mapping $\varphi$ is $L_p$-Lipschitz, and combining this with Assumption \ref{ass:sub_gauss} implies that by Theorem 1 in \citep{kontorovich2014} we have with probability at most $\exp(\frac{-\xi^2 n(\mathcal{T}_a)}{2\epsilon^2L_P^2\sigma^2})$ that the maximum difference between the empirical KL divergence and the true trajectory divergence is sufficiently far from its mean. Then using Lemma \ref{lem:rademach} we obtain an upper bound on this expected value with the appropriate constants. For a complete proof of the theorem please refer to Appendix \ref{ap:proof}. This theorem is useful because it indicates the empirical KL divergence derived from the MLE objective converges uniformly in probability to the true trajectory KL divergence.

\subsection{Concentration of Trajectory Divergence}
We can complete the proof of Theorem \ref{thm:concent_eq} using the results of Theorem \ref{thm:two_side_ineq} and the definition of the MLE. First, Theorem \ref{thm:two_side_ineq} implies that with high probability the trajectory divergence between the MLE parameters $\hat{\theta}_a,\hat{x}_{a,0}$ and true parameters $\theta^*_a,x_{a,0}^*$ is within $\mathcal{O}(\sqrt{\frac{d_x+d_\theta}{n(\mathcal{T}_a)}})$ of the empirical divergence between these two sets of parameters. Then, since $\hat{\theta}_a,\hat{x}_{a,0}$ minimize the empirical divergence and the empirical divergence of $\theta^*_a,x^*_{a,0}$ is zero, this means that the empirical divergence term is non-positive. Combining these two facts yields the concentration bound of Theorem \ref{thm:concent_eq}, and the complete proof is given in Appendix \ref{ap:proof}.

 
We conclude this section with an alternative statement of Theorem \ref{thm:concent_eq}.
%
\begin{corollary} \label{corr:alt_rep}
For $\alpha \in (0,1)$, with probability at least $1-\alpha$ we have
\begin{equation}\frac{1}{n(\mathcal{T}_a)}D_{a,\pi_1^{T}}(\theta^*_a,x^*_{a,0}||\hat{\theta}_a,\hat{x}_{a,0}) \leq B(\alpha)\sqrt{\frac{\log(1/\alpha)}{n(\mathcal{T}_a)}}.
\end{equation}
Where $B(\alpha) = \frac{c_f(d_x,d_\theta)}{\sqrt{\log(1/\alpha)}}+L_p\sigma\sqrt{2}.$
\end{corollary}
This result can be obtained by making the substitution $\xi = L_p\sigma\sqrt{\frac{\log(1/\alpha)}{n(T_a)}}$ into the expression in Theorem \ref{thm:concent_eq}. This corollary is significant because it allows us to derive confidence bounds for our parameter estimates with regards to their trajectory divergence. Note that the term $B(\alpha)$ differs from the term that would be derived by Chernoff-Hoeffding bounds applied to i.i.d. random variables by the addition of $\frac{c_f(d_x,d_\theta)}{\sqrt{\log(1/\alpha)}}$ to the standard variance term. The reason for this addition is that since we are using MLE for our parameter estimation our estimates will be biased, and this bias must be accounted for in the confidence bounds. Though there may exist specific models where MLE can provide unbiased estimates, we will only present analysis for the more general case.

\section{Policies for Optimizing ROGUE Models}
\label{sec:optim_pol}
This section develops two different policies which could be used to achieve effective regret bounds for ROGUE bandit models. These are a ROGUE Upper Confidence Bounds (ROGURE-UCB) policy, and an $\epsilon$-greedy policy ($\epsilon$-ROGUE). While both policies exhibit good theoretical performance, each policy has different advantages during implementation, because of how they trade off exploration-exploitation and computation. $\epsilon$-ROGUE is computationally efficient to implement, but during our experiments it tended to over-explore in the short term for many instances of ROGUE models. This may be because $\epsilon$-greedy policies explore uniformly, and even in the stationary case $\epsilon$-greedy policies are sensitive to their tuned parameters and do not perform well with multiple sub-optimal choices \citep{auer2002finite}.  On the other hand, ROGUE-UCB is a more computationally intensive policy; however, in our experiments it was capable of finding optimal arms with fewer samples then $\epsilon$-ROGUE. This makes ROGUE-UCB  more suitable for settings that require faster parameter identification such as those encountered in healthcare, while $\epsilon$-ROGUE is better suited for low risk but high volume settings such as online advertising.

Alhough several upper confidence bounds (UCB) policies have been proposed in the non-stationary setting \citep{garivier2008,besbes2014}, these existing policies provide regret of order $\mathcal{O}(\sqrt{T\log T})$ in the general adversarial case, but can achieve $\mathcal{O}(\log T)$ regret under certain model-based assumptions.  In contrast, both the ROGUE-UCB and $\epsilon$-ROGUE policies we construct achieve regret of order $\mathcal{O}(\log{T})$ under more general assumptions, which is optimal in the sense that it matches the lowest achievable rate for approximate policies in the stationary stochastic setting. We theoretically analyze the cumulative regret of both of these algorithms in the stochastic setting, and also develop a parameter free bound on regret for ROGUE-UCB. The algorithms we propose in this paper are designed for the stochastic ROGUE bandit setting.  We leave the discussion of the adversarial ROGUE bandit setting for future research.

\subsection{ROGUE Upper Confidence Bounds (ROGUE-UCB) Policy}
\label{sec:ducb}
%
 Pseudocode for ROGUE-UCB is given in Algorithm \ref{alg:s_ucb}, and the algorithm is written for the situation where the policy chooses actions over the course of $T$ time periods labeled $\{1,...,T\}$.  The upper confidence bounds used in this algorithm are computed using the concentration inequality from Theorem \ref{thm:concent_eq}. Much like other UCB policies, for the first $|\mathcal{A}|$ time steps of the algorithm each action $a$ will be tried once. Then after this initialization, at each time step, we will first compute the MLE estimates of the parameters for each action (i.e., $(\theta_a,\hat{x}_{0,a}) \forall a\in \mathcal{A}$) and then use Theorem \ref{thm:concent_eq} to form the upper confidence bound on the value of $g(\theta_a,x_{t,a})$, which we call $g^{UCB}_{a,t}$. Our approach for forming these bounds is similar to the method first proposed by \cite{garivier2011} for the KL-UCB algorithm used for stationary bandits. Here, since we know that with high probability the true parameters belong to $\mathcal{X}$ and $\Theta$, we find the largest possible value of $g(\theta_a,x_{t,a})$ within these sets. Finally, we choose the action that has the largest upper confidence bound, observe the result, and repeat the algorithm in the next time step. Methods for solving the component optimization problems effectively are highly dependent upon the underlying model and dynamics. In general when $h_a$ is nonlinear these problems will be non-convex; however, the MLE problem can be computed using dynamic programming or filtering approaches while the UCB problem can be computed directly or through a convex relaxation as a relaxation solution will still be a valid UCB. In implementation, if $d_x,d_\theta$ are small, then enumeration and fine grid search methods can be employed. We leave the discussion of more efficient algorithms for specific ROGUE models for future research.
%

The key theoretical result about the ROGUE-UCB algorithm concerns the regret $R_\Pi(T)$ of the policy computed by the ROGUE-UCB algorithm.
\begin{theorem} \label{thm:regret_them}
The expected regret $\mathbb{E}R_{\Pi}(T)$ for a policy $\Pi$ computed by the ROGUE-UCB algorithm is
\begin{equation}
\mathbb{E}R_\Pi(T) \leq L_g\diam(\mathcal{X}\times\Theta)\sum_{a\in\mathcal{A}} \Bigg(A(|\mathcal{A}|)^2 \frac{4\log T}{\delta_a^2} + \frac{\pi^2}{3}\Bigg).
\end{equation}
where $A(x) = B(x^{-4})$, and 
\begin{equation}
\begin{aligned}
\delta_a &= \min\{ \frac{1}{n(\mathcal{T}_a)}D_{a,\pi_1^T}(\theta_a,x_{a,0}||\theta_{a'},x_{a',0}): |g(h_a^t(x_{a,0}),\theta_a) - g(h_a^t(x_{a',0}),\theta_{a'})| \geq \frac{\epsilon_a}{2} \}\\
\epsilon_a &= \min_{a'\in\mathcal{A}\setminus a,t }\{|g(\theta_a,h_a^t(x_{a,0})) -g(\theta_{a'},h_a^t(x_{a',0}))|: g(\theta_a,h_a^t(x_{a,0})) \neq g(\theta_{a'},h_a^t(x_{a',0}))\}
\end{aligned}
\end{equation}
are finite and strictly positive constants.
\end{theorem}
\begin{remark}
This corresponds to a rate of order $\mathcal{O}(\log T)$ when $\liminf_T\delta_a > 0$.  In fact, $\liminf_T\delta_a > 0$ for many settings such as (with appropriate choice of model parameter values) the ROGUE GLM and ROGUE agent defined in Section \ref{sec:num_exp}.
\end{remark}
\begin{algorithm}
\begin{algorithmic}[1]
\caption{Reducing or Gaining Unknown Efficacy Upper Confidence Bounds (ROGUE-UCB)}\label{alg:s_ucb}
 \For{$t \leq |\mathcal{A}|$}
 \State $\pi_t = a$ such that $a$ hasn't been chosen before
\EndFor
\For{$|\mathcal{A}|\leq t \leq T$}
\For{$a \in \mathcal{A}$}
\State Compute: $\hat{\theta}_a,\hat{x}_{a,0} = \argmin\{ -\sum_{t\in \mathcal{T}_a}\log p(r_{a,t}|\theta_a,x_{a,t}): x_{a,t+1} = h_a(x_{a,t},\pi_{a,t}) \forall t \in 0,...,T \}$ 
\State Compute: $g^{UCB}_{a,t} = \max_{ \theta_a,x_{a,0}\in \Theta \times \mathcal{X}} \{g(\theta_a,h_a^t(x_{a,0})): \frac{1}{n(\mathcal{T}_a)}D_{a,\pi_1^{T}}(\theta_a,x_{a,0}||\hat{\theta}_a,\hat{x}_{a,0}) \leq A(t)\sqrt{\frac{4\log(t)}{n(\mathcal{T}_a)}}\}  $
\EndFor
\State Choose $\pi_t = \argmax_{a\in\mathcal{A}} g^{UCB}_{a,t}$
\EndFor
\end{algorithmic}
\end{algorithm}

To prove Theorem \ref{thm:regret_them}, we first present two propositions. The first proposition bounds the expected regret $R_\Pi(T)$ by the number of times an action is taken while it is suboptimal.
\begin{proposition} \label{prop:init_reg_bnd}
For a policy $\Pi$ calculated using the ROGUE-UCB algorithm, if $\tilde{T}_a = \sum_{t=1}^T \mathbf{1}\{\pi_t = a, a \neq \pi^*_t\}$, then $\mathbb{E}R_\Pi(T) \leq L_g\diam(\mathcal{X}\times\Theta)\sum_{a\in\mathcal{A}} \mathbb{E}\tilde{T}_a$.
\end{proposition}
For this proposition, we first use Assumption \ref{ass:lip_ass} to upper bound the value of the regret with respect to the $L_g$ and the diameter of the parameter set. Then since we are left with a finite sum of positive numbers, we can rearrange the summation term to obtain the expected number of suboptimal actions. For the detailed proof, please see Appendix \ref{ap:proof}. Next we proceed to prove a bound on the expected number of times a suboptimal action will be chosen.
\begin{proposition} \label{prop:wrong_pulls}
For a policy $\Pi$ calculated using the ROGUE-UCB algorithm, we have that $\mathbb{E} \tilde{T}_a \leq A(|\mathcal{A}|)^2 \frac{4\log T}{\delta_a^2} + \frac{\pi^2}{3}$, where $A(t) = B(t^{-4})$, $\delta_a = \min\{ \frac{1}{n(\mathcal{T}_a)}D_{a,\pi_1^T}(\theta_a,x_{a,0}||\theta_{a'},x_{a',0}): |g(h_a^t(x_{a,0}),\theta_a) - g(h_a^t(x_{a',0}),\theta_{a'})| \geq \frac{\epsilon_a}{2} \}$, and $\epsilon_a = \min_{a'\in\mathcal{A}\setminus a,t }\{|g(\theta_a,h_a^t(x_{a,0})) -g(\theta_a,h_a^t(x_{a,0}))|: g(\theta_a,h_a^t(x_{a,0})) \neq g(\theta_a,h_a^t(x_{a,0}))\}$.
\end{proposition}
To prove this proposition, we proceed in a manner similar to the structure first proposed by \cite{auer2002finite}. We must show that if an action is chosen at a time when it is suboptimal, then this implies that either we have not properly estimated its parameters (i.e., have not explored enough) or the true values of the parameters $x_{a,0},\theta_a$ or $x_{\pi^*_t,0},\theta_{\pi^*_t}$ are not contained inside their confidence bounds. Using these facts, we use Theorem \ref{thm:concent_eq} to show that the probability that all of these events occurring simultaneously is bounded, and then upper bound the expected number of times these events can occur.  Combining the results of Propositions \ref{prop:init_reg_bnd} and \ref{prop:wrong_pulls}, we thus prove the desired result of Theorem \ref{thm:regret_them}.  The full proofs of Proposition \ref{prop:wrong_pulls} and Theorem \ref{thm:regret_them} are provided in Appendix \ref{ap:proof}.

Developing a regret bound in the ROGUE setting that is free of problem-dependent parameters is quite challenging. This is because in contrast to the stationary setting, if $\epsilon_a,\delta_a$ are dependent on $T$ and can arbitrarily approach zero in the long run, the bound presented in Theorem \ref{thm:regret_them} might not hold. This behavior, however, depends greatly on the rate at which  $\epsilon_a,\delta_a$ approach zero. For instance,  it is clear that if $\epsilon_a = o(\frac{1}{T}), \forall a \in \mathcal{A}$ then the bound from Theorem \ref{thm:regret_them} might in fact be loose, so to achieve the worst case rate $\epsilon_a,\delta_a$ must approach zero quite slowly. In Section \ref{sec:eps_delt_discuss} we provide additional discussion of how these rates may impact regret and can be computed in certain settings.  In this section, we provide the following problem-dependent parameter free bound for the case when $\epsilon_a,\delta_a$ are bounded from below away from zero:

\begin{corollary} \label{cor:parameter_free}
If there exists $\epsilon$ such that $\forall a \in \mathcal{A}$ and $\forall T$ $\epsilon_a \geq \epsilon$, then in the worst case ROGUE-UCB achieves  $\mathbb{E}R_\Pi(T) = \mathcal{O}(T^{\frac{4}{5}}\big(|\mathcal{A}|A(|\mathcal{A}|)^2\sigma^2\log T\big)^{\frac{1}{5}})$.
\end{corollary}

The proof of this corollary can be found in Appendix \ref{ap:proof}, but here we present a sketch. Essentially, if the $\epsilon_a$ are bounded from below by $\epsilon$, then for small $\epsilon$ we would expect the regret to behave in the worst case linearly, namely at $\epsilon T$, and not according to the $\mathcal{O}(\log T)$ bound presented in Theorem \ref{thm:regret_them}. Thus, the result follows from finding an $\epsilon$ that minimizes both of these regret bounds and substituting into the given expressions. Unlike the stationary case, for which UCB1 has a worst case rate of $\mathcal{O}(\sqrt{T\log T})$, we note that this is a $\mathcal{O}((T^4\log T)^{\frac{1}{5}})$ worst case rate, and so is slightly worse but still sub-linear. This is not surprising, and can be attributed to the high variability of the problem parameters. In fact, this agrees with the bounds computed by \cite{besbes2014}, who show that if total variation is proportional to $T$, we expect worst case regret of the form $\mathcal{O}(T^{\frac{2 +\beta}{3}}) $ for some $\beta \in (0,1]$, a rate our result meets up to a $\log$ factor.

\subsection{$\epsilon$-Greedy Policy for ROGUE Bandits}
\label{sec:eps_greed}
In this section we develop an $\epsilon$-greedy policy that can be used to optimize ROGUE bandits. The pseudo code of this method can be found in Algorithm \ref{alg:eps_greedy}. At each time step $t$, the algorithm samples a standard uniform random variable $u_t$. If $u_t$ is above some threshold $\epsilon_t$ the algorithm performs a pure exploration step (i.e. chooses an arm from $\mathcal{A}$ uniformly randomly). Otherwise, the algorithm performs a greedy optimization step by computing the MLE values  for the expected rewards of each of the arms, and then playing the arm with the highest MLE expected reward. Note that $\epsilon_t$ decreases as the bandit is played -- this is critical to ensure that as the MLE estimates improve the algorithm makes fewer unnecessary explorations. Additionally, if $\epsilon_t$ is held at some constant value $\epsilon$, it is easy to see that this results in linear time regret, specifically a lower bound regret of $\Omega(\epsilon T)$.

We next show the following theoretical results for the cumulative regret performance of $\epsilon$-ROGUE:

\begin{algorithm}
	\begin{algorithmic}[1]
		\caption{Reducing or Gaining Unknown Efficacy $\epsilon$-Greedy ($\epsilon$-ROGUE)}\label{alg:eps_greedy}
		\State Set $c,d$
		\For{$1 \leq t \leq T$}
		\State Set: $\epsilon_t = \min\{1, \frac{c|\mathcal{A}|}{\delta_{min}^2t}\}$
		\State Sample: $u_t \sim \mathcal{U}(0,1)$
		\If {$\epsilon_t \leq u_t$}
			\State  Sample: $\pi_t \sim \text{Cat}(\mathcal{A}) $
			\State \textbf{Continue}
		\Else
			\For{$a \in \mathcal{A}$}
				\State Compute: $\hat{\theta}_a,\hat{x}_{a,0} = \argmin\{ -\sum_{t\in \mathcal{T}_a}\log p(r_{a,t}|\theta_a,x_{a,t}): x_{a,t+1} = h_a(x_{a,t},\pi_{a,t}) \forall t \in 0,...,T \}$ 
			\State Compute: $\hat{g}_{a,t} = g(\hat{\theta}_a,h_a^t(\hat{x}_{a,0}))$
			\EndFor
			\State Choose $\pi_t = \argmax_{a\in\mathcal{A}} \hat{g}_{a,t}$
		\EndIf
		\EndFor
	\end{algorithmic}
\end{algorithm}

	\begin{theorem} \label{thm:eps_greed_reg}
		For any time $T> \frac{c|\mathcal{A}|}{\delta^2_{min}} $, the probability a sub-optimal arm is pulled for a policy $\Pi$ computed by the $\epsilon$-ROGUE algorithm is at most:
		\begin{equation}
		\frac{8L_p^2\sigma^2}{\delta^2_{min}}\exp\big(\frac{c_f(d_x,d_\theta)^2}{2L_p^2\sigma^2}\big)\big(\frac{eT\delta_{min}^2}{|\mathcal{A}|c}\big)^\frac{-c}{8L_p^2\sigma^2} + \frac{c}{\delta_{min}^2}(\log \frac{eT\delta_{min}^2}{|\mathcal{A}|c})\big(\frac{eT\delta_{min}^2}{|\mathcal{A}|c}\big)^\frac{-3c}{28\delta_{min}^2} + \frac{c}{\delta_{min}^2T}
		\end{equation}
		Where $\delta_{min} \leq \min_{a \in \mathcal{A}} \delta_a$.
	\end{theorem}

	\begin{remark}
		This instantaneous regret bound implies asymptotic cumulative regret of order $\mathcal{O}(\log(T))$ since the first several terms are $o{\tiny }(\frac{1}{T})$ for a sufficiently large choice of constant $c$, namely $c \geq \max\{8L_p^2\sigma^2, \frac{28\delta_{min}^2}{3}\}$. 
	\end{remark}
	
	A complete proof for Theorem \ref{thm:eps_greed_reg} can be found in Appendix \ref{ap:proof}, but here we provide a sketch of the proof. Essentially, to bound the probability of a sub optimal arm being pulled, we first use the union bound to bound it from above by the probability it was pulled randomly in an exploration step or during an exploitation step.  We then bound the probability of subotimality in the exploration step by using the concentration bounds from Theorem \ref{thm:concent_eq}.

\subsection{Discussion of the Behavior of $\epsilon_a$ and $\delta_a$}
\label{sec:eps_delt_discuss}
The constants $\epsilon_a$ and $\delta_a$ are critical for the results presented in the previous two sections. However, these quantities have a non-trivial dependence on the underlying model of the ROGUE bandit as well as the algorithm utilized, and are thus difficult to analyze in the most general case. Below we present two different examples meant to illustrate how $\epsilon_a$ and $\delta_a$ may behave for certain models and how this may impact performance.

\subsubsection{Exponentially Small $\epsilon_a,\delta_a$} For this section, consider a two-armed ROGUE bandit model with arms $\mathcal{A}:= \{0,1\}$ such that $\theta_1=\theta_0=0$ and $r_{a,t} \sim \mathcal{N}(x_{a,t},1) \forall a \in \mathcal{A}, \mathcal{X}:=[0,1]$, hence $g(\theta_a,x_{a,t}) = x_{a,t}$. Define the dynamics for each arm as follows:
\begin{equation}
\begin{aligned}
&  h_0(x_{0,t},\pi_{0,t}) = \text{proj}_{[0,1]}\big( 0 \cdot x_{0,t} + 0 \cdot \pi_{0,t}\big) \\
&h_1(x_{1,t},\pi_{1,t}) = \text{proj}_{[0,1]}\big( 0.5 \cdot x_{1,t} + 0 \cdot \pi_{1,t}\big) 
\end{aligned}
\end{equation}

Observe that here arm 0 stays at state $x_{0,t} = 0$ for all time periods $t$, while the state of arm one can be described by $x_{1,t} = \big(\frac{1}{2}\big)^t x_{1,0}$. So if $x_{1,0} > 0$.  As time progresses, the expected reward of each arm get exponentially close, and hence $\epsilon_a$ decreases exponentially. Moreover, using the KL divergence of normal distributions, it is clear that  $\delta_a = \Omega(\big(\frac{1}{2}\big)^{2T} x_{1,0}^2)$. This means that asymptotically the regret bounds provided in Sections \ref{sec:ducb} and \ref{sec:eps_greed} would no longer be of order $\mathcal{O}(\log T)$.


\subsubsection{Constant Bounded $\epsilon_a,\delta_a$}
As before, consider a two armed ROGUE bandit model with arms $\mathcal{A} := \{0,1\}$ such that $\theta_1 = \theta_0 = 0$ and $r_{a,t} \sim \mathcal{N}(x_{a,t},1)$ $\forall a \in \mathcal{A}$, $\mathcal{X} := [0,1]$, and $g(\theta_{a,t},x_{a,t}) =x_{a,t}$. Furthermore define the dynamics of each arm as follows:
\begin{equation}
\begin{aligned}
& h_0(x_{0,t},\pi_{0,t}) = \text{proj}_{[0,1]}\big(x_{0,t} + 2 \cdot \pi_{0,t} - 1\big) \\
& h_1(x_{1,t},\pi_{1,t}) = \text{proj}_{[0,1]}\big(x_{1,t} - 2 \cdot \pi_{0,t} + 1\big) 
\end{aligned}
\end{equation}

Here, the dynamics are only with respect to the action taken on arm 0 since in a two armed bandit setting, by definition $\pi_{0,t} = 1 - \pi_{1,t}$.  From these dynamics it is clear that regardless of which arm is played, whenever an action is taken the reward of of the played arm will decrease to 0 while the reward of the unplayed arm will increase to 1. Hence, for any horizon $T$, in this case $\epsilon_0 = \epsilon_1 = \min\{1, |x_{0,0} - x_{1,0}|\}$ if $x_{0,0} \neq x_{1,0}$ and $\epsilon_0 = \epsilon_1 = 1$ otherwise. Using the KL divergence between normal distributions with different means and equal variances, this implies that $\delta_0 = \delta_1 \geq \frac{1}{8}\min\{1, (x_{0,0} - x_{1,0})^2 \}$ when $x_{0,0} \neq x_{1,0}$ and $\delta_0 = \delta_1 = \frac{1}{8}$ otherwise. for any horizon $T$. This holds asymptotically, and so since $\delta_a,\epsilon_a$ are bounded asymptotically, this indicates that the regret results from Sections \ref{sec:ducb} and \ref{sec:eps_greed} will hold and regret will be of order $\mathcal{O}(\log T)$.

%
\section{Numerical Experiments}
\label{sec:num_exp}
In this section, we perform two numerical experiments where the policies computed by the ROGUE-UCB and $\epsilon$-ROGUE algorithms are compared against the policies computed by other non-stationary bandit algorithms. The first experiment considers the ROGUE GLM described in Section \ref{sec:dglm}, and  specifically looks at the logistic regression instantiation of ROGUE GLM. We use synthetically generated data for this first experiment. Next, we perform an experiment in the context of healthcare-adherence improving interventions to increase physical activity, which can be modeled using the ROGUE agent from Section \ref{sec:stack_band}. Using real world data from the mDPP trial \citep{fukuoka2014b}, we show how ROGUE-UCB and $\epsilon$-ROGUE can be implemented to personalize messages for participants in this intervention. All experiments in this section were run using Python 3.5.2 and Anaconda on a laptop computer with a 2.4GHz processor and 16GB RAM.

Our goal with these experiments is to show how each of the bandit algorithms performs in terms of long run reward and how quickly the algorithms can learn the underlying system state. While these problems are equivalent in the stationary bandit setting, in the non-stationary setting incorrect actions taken due to improper system identification may still have rewards that are quite close to optimal. As such, we consider two different performance measures for each algorithm, the tracking regret which measures how quickly the unknown model parameters are identified, and the long run reward. Our results show that both the Tuned ROGUE-UCB and $\epsilon$-ROGUE algorithms outperform the other candidate algorithms both in terms of long run average reward and average regret. However, while both ROGUE algorithms have similar performance in the long run, ROGUE-UCB tends to identify model parameters faster then $\epsilon$-ROGUE as it performs less exploration but at a slightly greater computational cost at each time step. This suggests that ROGUE-UCB is most useful for settings that require fast identification, such as settings that require higher consideration for risk and fairness (e.g., healthcare settings), while $\epsilon$-ROGUE is ideal for low risk settings where computational resources are the main implementation constraint (e.g., targeted advertising). 

\subsection{Tuned ROGUE-UCB}
As has been noted for other UCB policies \citep{auer2002finite,garivier2008,bouneffouf2016}, the high probability bounds derived theoretically for these methods are often too conservative. While the $\mathcal{O}(\sqrt{\frac{\log{t}}{n(\mathcal{T}_a)}})$ is a tight rate, the term $A(t)$ is too conservative. Drawing inspiration from \cite{auer2002finite} who used asymptotic bounds for Tuned UCB, we similarly construct a variant of our algorithm: This variant is described in Algorithm \ref{alg:tuned_ducb} and called Tuned ROGUE-UCB. Using the results of \cite{shapiro1993asymptotic}, we note that if the MLE $\hat{\theta}_a,\hat{x}_{a,0}$ are in the interior of the feasible region and are consistent, then they are asymptotically normally distributed with a variance equal to their Fisher information. Using these results and the delta method \citep{qu2011}, we can derive the quantity $\mathcal{S}_{a,\pi_1^T}(\theta_a,x_{a,0}||\hat{\theta}_a,\hat{x}_{a,0}) = \frac{1}{n(\mathcal{T}_a)^2}\nabla_{\theta',x'}D_{a,\pi_1^{T}}(\theta_a,x_{a,0}||\theta',x')^T \mathcal{I}_{\{r_t\}_{t\in\mathcal{T}_a}}(\theta',x')^{-1}\nabla_{\theta',x'}D_{a,\pi_1^{T}}(\theta_a,x_{a,0}||\theta',x')\big|_{\theta',x'=\hat{\theta}_a,\hat{x}_{a,0}}$, which is the asymptotic variance of the average trajectory KL-Divergence.  Here, $\eta$ is a constant that corresponds to the maximum value of the KL-divergence; $\mathcal{I}_{\{r_t\}_{t\in\mathcal{T}_a}}(\theta',x')$ represents the observed trajectory Fisher information, which can be calculated as $\mathcal{I}_{\{r_t\}_{t\in\mathcal{T}_a}}(\theta',x') = \sum_{t\in\mathcal{T}_a}\mathcal{I}_{r_t}(\theta',x')$, due to Assumption \ref{ass:indep_ass}. As an implementation note, if the empirical information matrix is singular, then the Moore-Penrose pseudoinverse should be used to achieve similar asymptotic results \citep{hero1997}. Note that although these asymptotic bounds work well in practice, they are not high probability bounds and do not provide the same theoretical guarantees as the ROGUE-UCB algorithm. A full analysis of regret for Tuned ROGUE-UCB is beyond the scope of this work. Instead, we only consider empirical analysis of this algorithm to show its strong performance. 

\begin{algorithm}[h]
\begin{algorithmic}[1]
\caption{Tuned ROGUE-UCB}\label{alg:tuned_ducb}
 \For {$t \leq |\mathcal{A}|$}
 	\State $\pi_t = a$ such that $a$ hasn't been chosen before
\EndFor
\For {$|\mathcal{A}|\leq t \leq T$}
	\For {$a \in \mathcal{A}$}
		\State Compute: $\hat{\theta}_a,\hat{x}_{a,0} = \argmin\{ -\sum_{t\in \mathcal{T}_a}\log p(r_{a,t}|\theta_a,x_{a,t}): x_{a,t+1} = h_a(x_{a,t},\pi_{a,t}) \forall t \in 0,...,T \}$ 
		\State Compute: $g^{UCB}_{a,t} = \max_{ \theta_a,x_{a,0}\in \Theta \times \mathcal{X}} \Big\{g(\theta_a,h_a^t(x_{a,0})): \frac{1}{n(\mathcal{T}_a)}D_{a,\pi_1^{T}}(\theta_a,x_{a,0}||\hat{\theta}_a,\hat{x}_{a,0}) \leq \sqrt{\min\{\frac{\eta}{4},\mathcal{S}_{a,\pi_1^T}(\theta_a,x_{a,0}||\hat{\theta}_a,\hat{x}_{a,0})\}\frac{\log(t)}{n(\mathcal{T}_a)}}\Big\}$
	\EndFor
\State Choose $\pi_t = \argmax_{a\in\mathcal{A}} g^{UCB}_{a,t}$
\EndFor
\end{algorithmic}
\end{algorithm}

\subsection{Experimental Design}
We examined two settings for our experiments, which correspond to the instantiations of ROGUE bandits presented in Sections \ref{sec:stack_band} and \ref{sec:dglm}. For each of the scenarios, we compared the Tuned ROGUE-UCB algorithm and $\epsilon$-ROGUE algorithms with $\epsilon_t = \frac{1}{t}$ to policies determined by six alternative methods. For each scenario, we present two result metrics: cumulative regret of each algorithm in that scenario and the average reward to date of the algorithm. While these two measures are related, a key difference is that in the non-stationary setting sub-optimal actions may not have a significantly lower expected reward than the optimal action at all time periods. Hence, while an algorithm may incur a significant amount of regret it could still achieve a high amount of reward. The five alternative algorithms we used for comparison are as follows: 
\begin{enumerate}
\item {\bf Pure Exploration:} First, we considered a completely random, or ``pure exploration'' algorithm, which chooses an action uniformly at random from the set of available actions. 

\item {\bf Stationary Upper Confidence Bound (UCB1):} Next, we considered the UCB1 algorithm \citep{auer2002finite}, which is designed for stationary bandits. This approach uses the sample average as an estimate of the expected reward of each action and utilizes a padding upper confidence bound term derived from Hoeffding's bound. In our experiments, we implemented Tuned UCB1 \citep{auer2002finite}, which replaces the theoretical constants by the asymptotic variance of the sample average and a small constant that corresponds to the maximum variance of a Bernoulli random variable (since the rewards are bounded between 0 and 1).

\item {\bf Discounted Upper Confidence Bounds (D-UCB):} D-UCB is an upper confidence bound approach designed for non-stationary systems. It utilizes an exponentially weighted average of the reward observations to estimate the expected reward at the current time period and a square root padding function to provide upper confidence bounds \citep{garivier2008}. The weighted average is constructed with a positive discount factor that decreases the influence of older observations on the reward estimate to zero as time goes on. We implemented this algorithm with its optimal theoretical parameters, as described in \cite{garivier2008}.
 
\item {\bf Sliding Window Upper Confidence Bounds (SW-UCB):} The next approach we considered is the SW-UCB approach. This algorithm considers a fixed window size of how many action choices to ``keep in memory'', and computes the estimate of the expected action rewards as the average of these choices \citep{garivier2008}. We implemented this algorithm with its optimal theoretical parameters as proposed by \cite{garivier2008}.

\item {\bf Exploration and Exploitation with Exponential Weights (EXP3):} The next bandit algorithm we considered in our experiments is the EXP3 algorithm. Essentially, EXP3 is a modification of the exponential weights algorithm used in online optimization to the bandit setting where not all action rewards are observed \citep{auer2002}. Though EXP3 is designed for stationary bandits, unlike UCB approaches that assume a stochastic setting, it is meant for  adversarial bandits, which makes it potentially robust to non-stationarity. The particular variant of EXP3 we utilized is EXP3.S proposed by \cite{auer2002}, which is designed for arbitrary reward sequences, using the theoretically optimal parameters as proposed by the authors. Specifically, the learning rate $\gamma$ was set to $\min\Big\{1,\sqrt{\frac{|\mathcal{A}|}{T}\log(|\mathcal{A}|T)} \Big\} $ and prior weighting $\alpha$ was set to $\frac{1}{T}$. Essentially, this parameter tuning assumes the full horizon the bandit is to be run for is known. 

\item {\bf Restarting EXP3 (REXP3):} The last bandit algorithm we considered in our experiments is the REXP3 algorithm proposed by \cite{besbes2014}. This algorithm is a modification of the EXP3 algorithm for a non-stationary setting, where the none-stationarity has bounded total variation. This approach  performs EXP3 update steps for a preset number of arm pulls, and once this number is reached, the algorithm resets as if it is restarting. The length of the period before the reset depends on the total time horizon of the bandit, as well as the total variation of the underlying reward. In our experiments, we tuned this algorithm using various values for the total variation and found that the best results overall were achieved by setting this value to be proportional to the total horizon $T$, so we only present these results.

\end{enumerate}

\subsection{ROGUE Logistic Regression}
For this experiment, we consider the logistic regression instantiation of the ROGUE GLM presented in Section \ref{sec:dglm}.  Our setup includes two actions whose rewards $r_{t,a}$ are Bernoulli with a logistic link function of the form $g(x,\theta) = \frac{1}{1+ \exp(-a\theta - b x)}$. The initial parameters and dynamics matrices for each of the actions are presented in Table \ref{table:glm_params}. Here, the sets $\mathcal{X}$ and $\Theta$ were set to $[0,1]$. Action 0 has a significant portion of its reward dependent on the time varying state $x_t$, and recovers its reward slowly but also decreases slowly. On the other hand, Action 1 has more of its expectation dependent on the stationary component $\theta$, but it expectation decreases faster than that of Action 0. 
  
\begin{table}[h]
\begin{center}
\begin{tabular}{|c|c|c|c|c|c|c|c|}
\hline
Action & $x_0$ &$\theta$ & $A$ & $B$ & $K$ & $\alpha$ & $\beta$\\ \hline 
0   &  0.1  &   0.5 & 0.6 & -1.0 &  0.5 & 0.4 & 0.6 \\ \hline
1   &  0.3  &   0.7 & 0.7 & -1.2 &  0.5 & 0.7 & 0.3 \\ \hline
\end{tabular}
\end{center}
\caption{Experimental parameters for each action for the logistic ROGUE GLM simulation}
\label{table:glm_params}
\end{table}

\begin{figure}
\centering
\includegraphics[scale=0.9]{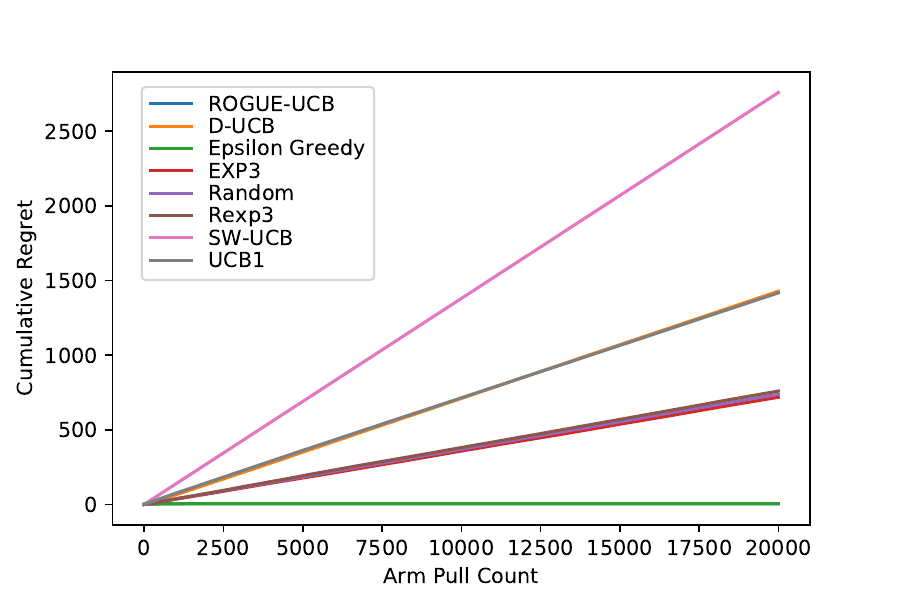}
\caption{Comparison of cumulative regret between the different bandit algorithms for the logistic ROGUE GLM.}
\label{fig:logit_regret}
\end{figure}

\begin{figure}
\centering
\includegraphics[scale=0.9]{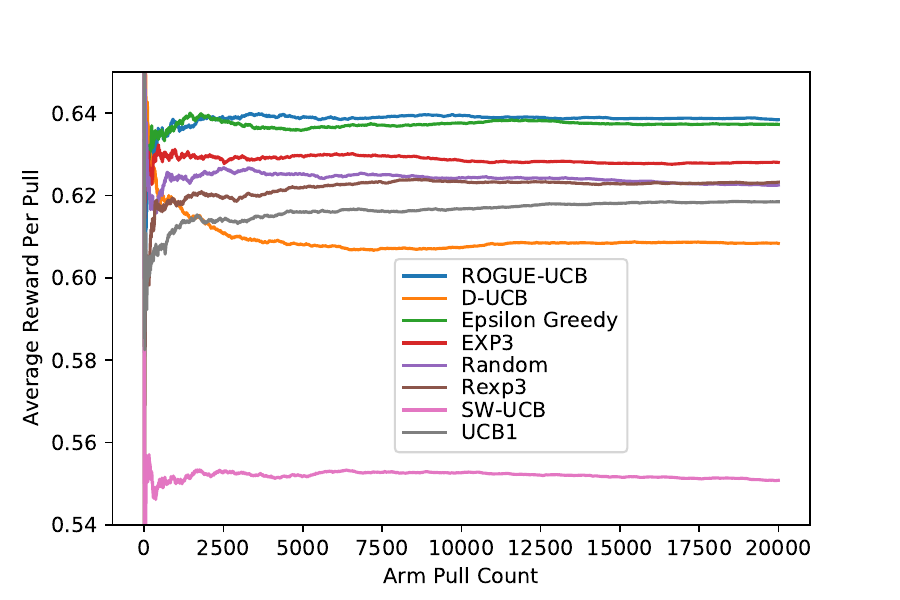}
\caption{Comparison of average reward to date between the different bandit algorithms for the logistic ROGUE GLM.}
\label{fig:logit_reword}
\end{figure}

The experiments were run for 20,000 action choices and replicated 30 times for each of the candidate algorithms. Figure \ref{fig:logit_regret} shows the cumulative regret accrued by each of the algorithms averaged across the replicates, and Figure \ref{fig:logit_reword} shows the average reward per action for each algorithm averaged across the replicates. As expected in these experiments, the UCB1 algorithm achieves linear regret since it assumes a stationary model and thus converges to a single action, which causes a large gap between the expectations of the two actions. Interestingly, SW-UCB and D-UCB also perform worse than random choices. A key note here is that D-UCB and SW-UCB assume that action rewards do not change frequently and are independent of the choices. However, D-UCB outperforms SW-UCB since the weighted average contains more information about the trajectory of the expected reward of each action while data from earlier choices are removed from the estimates in the sliding window. EXP3, REXP3, and random action selection perform approximately the same in terms of both regret and expected reward. This is unsurprising because the weighting scheme in EXP3 emphasizes the rewards of the past action states as opposed to current action states, and high total variation essentially reduces the restart horizon of REXP3 which means it retains less information about the overall reward trajectory. In terms of both regret and reward, Tuned ROGUE-UCB and $\epsilon$-Rogue substantially outperform the other approaches.

\begin{figure}
	\centering
	\includegraphics[scale=0.9]{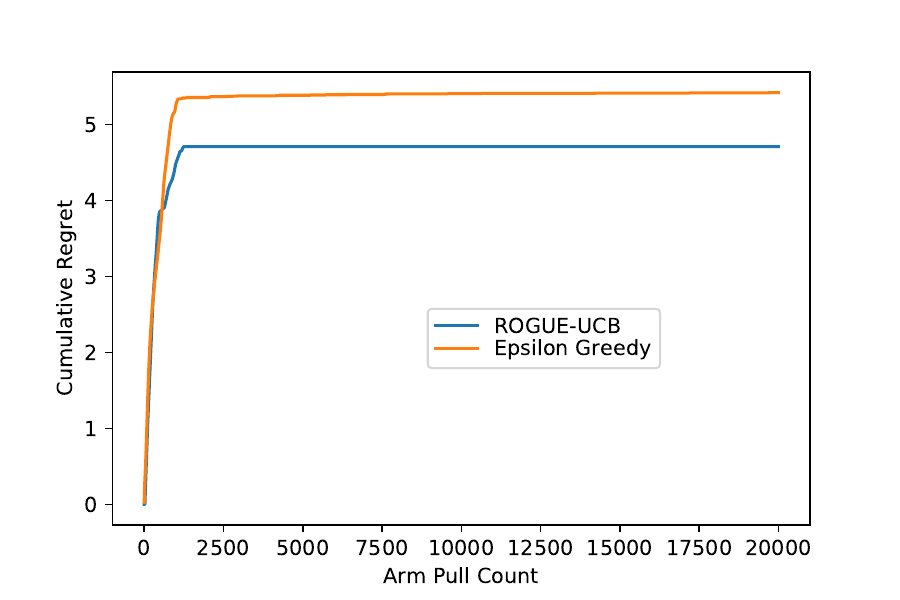}
	\caption{Comparison of cumulative regret between the different bandit algorithms for the logistic ROGUE GLM.}
	\label{fig:logit_regret_model}
\end{figure}

\begin{figure}
	\centering
	\includegraphics[scale=0.9]{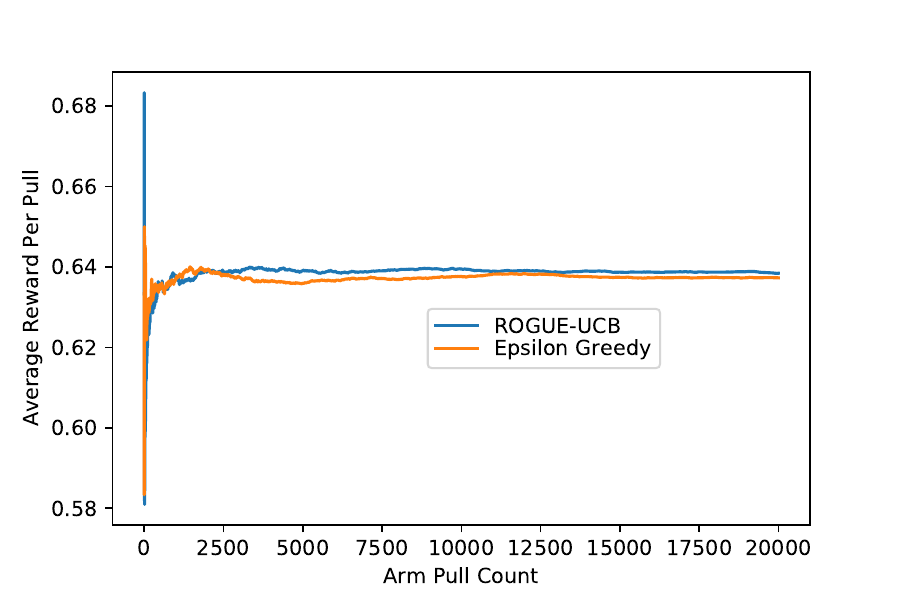}
	\caption{Comparison of average reward to date between the different bandit algorithms for the logistic ROGUE GLM.}
	\label{fig:logit_reword_model}
\end{figure}

In Figures \ref{fig:logit_regret_model} and \ref{fig:logit_reword_model} we compare the two model-based approaches in terms of their cumulative regret and average reward, respectively. While the other approaches seem to obtain linear regret, both ROGUE-UCB and $\epsilon$-ROGUE do in fact have regret on the order of $\mathcal{O}(\log T)$ in this experiment. Moreover, in this particular non-stationary setting we observe that these two algorithms have very similar long run average rewards while ROGUE-UCB obtains slightly better cumulative regret. This can be attributed to ROGUE-UCB performing less exploration then $\epsilon$-ROGUE.

\subsection{Healthcare-Adherence Improving Intervention for Increasing Physical Activity}
Next, we consider an experiment using real world data from the mobile diabetes prevention program (mDPP) \citep{fukuoka2014b}. This was a randomized control trial (RCT) that was conducted to evaluate the efficacy of a 5 month mobile phone based weight loss program among overweight and obese adults at risk for developing type 2 diabetes and was adapted from the Diabetes Prevention Program (DPP) \citep{dpp2002,dpp2009}. Sixty one overweight/obese adults were randomized into an active control group that only received an accelerometer (n=31) or a treatment group that received the mDPP mobile app plus the accelerometer and clinical office visits (n=30). Changes in primary and secondary outcomes for the trial were clinically and statistically significant. The treatment group lost an average of 6.2 $\pm$ 5.9 kg (-6.8\% $\pm$ 5.7\%) between baseline and the 5 month follow up while the control group gained 0.3 $\pm$ 3.0 kg (0.3\% $\pm$ 5.7 \%) (p $<$ 0.001). The treatment group's steps per day increased by 2551 $\pm$ 4712 compared to the control group's decrease of 734 $\pm$ 3308 steps per day (p $<$ 0.001). Additional details on demographics and other treatment parameters are available in \citep{fukuoka2014b}.

One key feature of the mDPP application was the ability for the clinicians to send daily messages to the participants to encourage that they adhere to the intervention and maintain a sufficiently increased activity level. Broadly speaking, there were 5 different message categories that the clinicians could choose to send to the patients. These categories are self-efficacy/confidence, motivation/belief/attitude, knowledge, behavior reinforcement, and social support. Each day the experimental group would receive a preprogrammed message from one of these categories, and all participants received the same messages each day. For our simulations, we used the data of what messages were sent to what participants, as well as their daily step counts.

\subsubsection{Patient Model}
For our experiment, we used a behavioral analytics model of patient behavior first proposed by \cite{aswani2016}. Here, each patient is assumed to be a utility maximizing agent who chooses how many steps to take each day based on previous behavior and the intervention implemented. We defined each of the different message categories be one of the actions of the bandit, which forms a ROGUE agent model as described in Section \ref{sec:stack_band}. Using the notation of Section \ref{sec:stack_band}, let $c_t$ be a sequence of i.i.d. Laplace random variables with mean zero and shape parameter $\theta$. This means $\sigma^2(\theta) = 2\theta^2$. After normalizing the step counts to be in $[0,1]$ (where 1 is equal 14,000 steps), we can then write the reward distribution of a particular message type $a$ as $p(r_t|\{x_{a,t},\theta_a\}) = \frac{1}{2}\exp(\frac{-x_{a,t}}{\theta_a}) \delta(r_t) + \frac{1}{2}\exp(\frac{x_{a,t}-1}{\theta_a}) \delta(1-r_t) + \frac{1}{2\theta_a} \exp(\frac{-|r_t-x_t|}{\theta_a})\mathbf{1}[r_t \in (0,1)]$, where the state $x_{a,t} \in [0,1]$ and $\theta_a \in [\epsilon,1]$ for a small $\epsilon >0$. This results in a reward function $g(x,\theta) = x + \frac{\theta}{2}(\exp(\frac{-x}{\theta})- \exp(\frac{x-1}{\theta}))$. Using Laplace noise has the advantage of allowing commercial mixed integer programming solvers to be used for offline parameter estimation by solving inverse optimization problems \citep{aswani2015,aswani2017,mintz2017}. Using this MILP reformulation and behavioral models, we estimated the respective trajectory parameters for each message group and each patient of the treatment group for which we had data. These initial parameters were found using the Gurobi Solver in Python \citep{gurobi}.

\subsubsection{Simulation Results}
\begin{figure}
\centering
\includegraphics[scale=0.9]{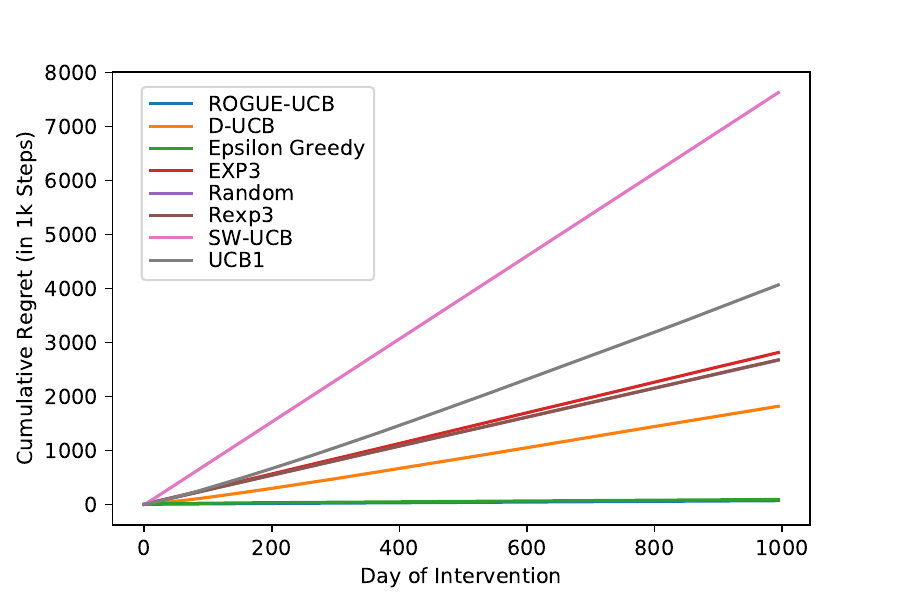}
\caption{Comparison of cumulative regret between the different bandit algorithms for the healthcare-adherence improving intervention.}
\label{fig:mdpp_regret}
\end{figure}

\begin{figure}
\centering
\includegraphics[scale=0.9]{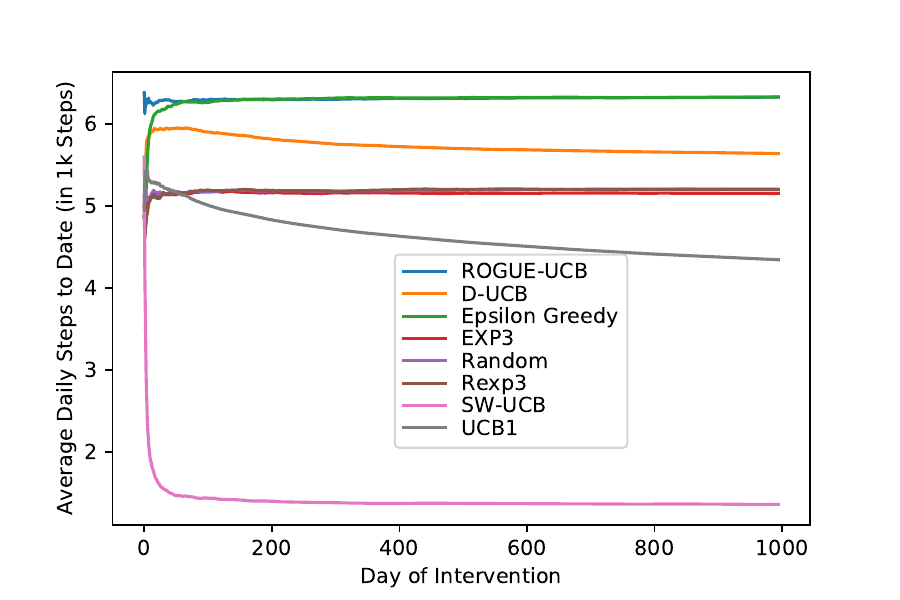}
\caption{Comparison of average reward to date between the different bandit algorithms for the healthcare-adherence improving intervention.}
\label{fig:mdpp_reword}
\end{figure}

This simulation was conducted using the mDPP data described above. Each experiment consisted of 1,000 action choices, which would correspond to about two years of a message based physical activity intervention, and 10 replicates of the simulation were conducted per patient and algorithm. The results in Figures \ref{fig:mdpp_regret} and \ref{fig:mdpp_reword} represent averages across all patients and replicates. Since we are using real data, the interpretation of the y-axis of each of the plots corresponds to number of steps in units of 1,000 steps, and the x-axis corresponds to the day of the intervention.

ROGUE-UCB and $\epsilon$-ROGUE outperform all other algorithms both in terms of regret and average reward. In terms of regret, both ROGUE-UCB and $\epsilon$-ROGUE obtain logarithmic-in-time regret; however, as before, ROGUE-UCB achieves lower cumulative regret than $\epsilon$-ROGUE. While D-UCB is the only other algorithm that can outperform pure exploration, it only obtains linear regret. In terms of average reward, ROGUE-UCB, $\epsilon$-ROGUE, and D-UCB are the only algorithms that outperform pure exploration. Interpreting these results in the healthcare context of this intervention, we find that the improved predictive model and use of MLE estimates within our specialized ROGUE Bandit algorithms results in an increase of 1,000 steps a day (approximately a half-mile more of walking per day) relative to the next best algorithm, which is a significant increase in activity. We note that in terms of average reward, both ROGUE-UCB and $\epsilon$-ROGUE have similar performance; however, ROGUE-UCB has higher initial rewards due to fewer exploration steps. This property is of particular interest in this healthcare setting, as it insures that each participating individual receives the most appropriate recommendations at a faster rate.

\begin{figure}
	\centering
	\includegraphics[scale=0.9]{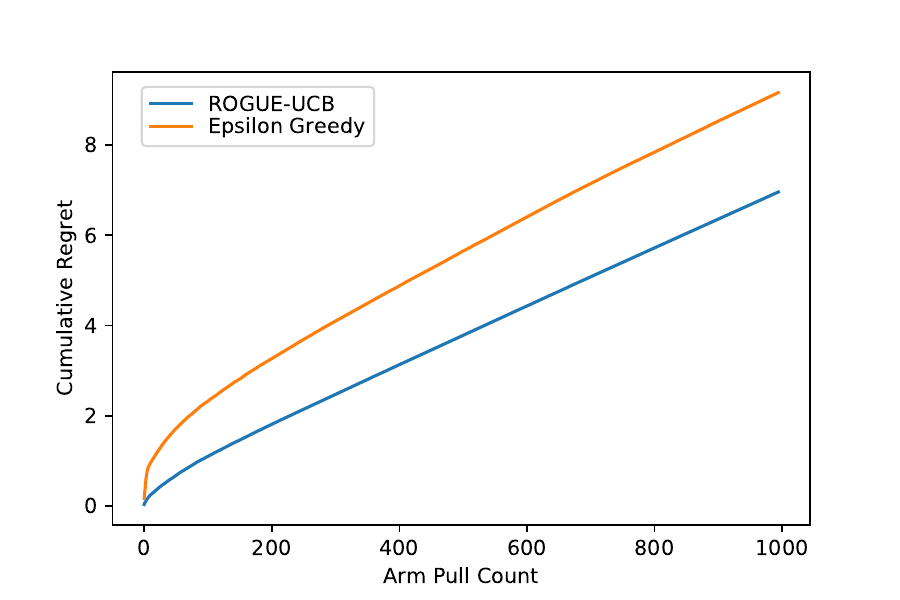}
	\caption{Comparison of cumulative regret between the different bandit algorithms for the healthcare-adherence improving intervention.}
	\label{fig:mdpp_regret_model}
\end{figure}
\begin{figure}
	\centering
	\includegraphics[scale=0.9]{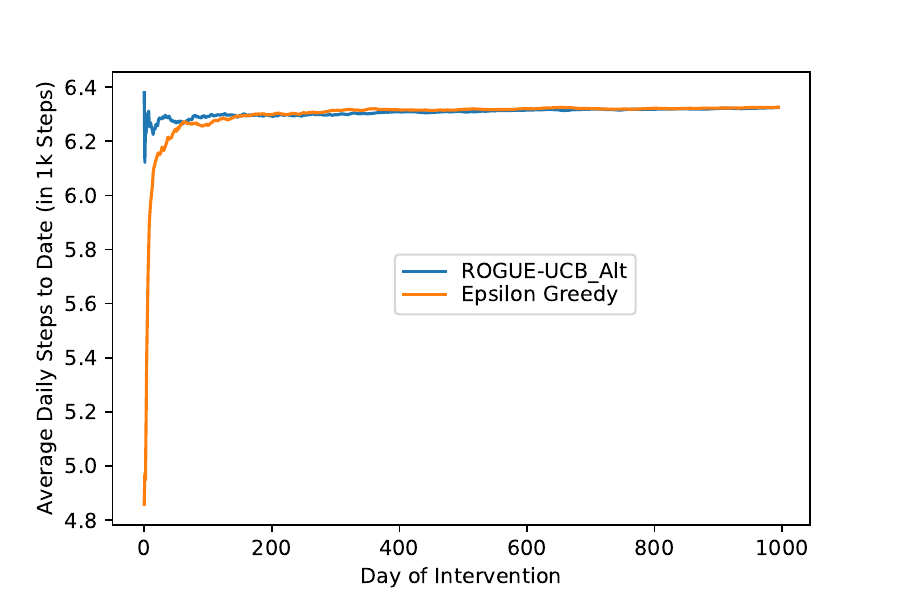}
	\caption{Comparison of average reward to date between the different bandit algorithms for the healthcare-adherence improving intervention.}
	\label{fig:mdpp_reword_model}
\end{figure}
\section{Conclusion}
In this paper, we defined a new class of non-stationary bandit models where the specific actions chosen influence the reward distributions of each action in subsequent time periods through a specific model. We conducted a finite sample analysis of the MLE estimates in this setting, and derived finite time concentration bounds. Using this analysis, we developed both the ROGUE-UCB  and $\epsilon$-ROGUE algorithms that provide policies for these bandit models. Our theoretical results show that in expectation ROGUE-UCB and $\epsilon$-ROGUE achieve logarithmic  in time regret. This is a substantial improvement over model-free algorithms, which can only achieve a square-root regret in the worst case scenario. We then showed through simulations using real and artificial data, that with minor modification, the ROGUE-UCB and $\epsilon$-ROGUE algorithms significantly outperform state of the art bandit algorithms both in terms of cumulative regret and average reward. These results suggest that ROGUE bandits have strong potential for personalizing health care interventions, and in particular for healthcare-adherence improving interventions.

\ACKNOWLEDGMENT{The authors gratefully acknowledge the support of NSF Award CMMI-1450963, UCSF Diabetes Family Fund for Innovative Patient Care-Education and Scientific Discovery Award, K23 Award (NR011454), and the UCSF Clinical and Translational Science Institute (CTSI) as part of the Clinical and Translational Science Award program funded by NIH UL1 TR000004.}

\bibliographystyle{informs2014} 
\bibliography{rogue}

\hrulefill
\vspace{1.5mm}

\textbf{Yonatan Mintz:} Yonatan is a Postdoctoral Research Fellow at the H. Milton Stewart School of Industrial and Systems Engineering at the Georgia Institute of Technology, previously he completed his PhD at the department of Industrial Engineering and Operations research at the University of California, Berkeley. His research interests focus on the application of machine learning and optimization methodology for personalized healthcare and fair and accountable decision making.

\textbf{Anil Aswani:} Anil is an assistant professor in the Department of Industrial Engineering and Operations Research at the University of California, Berkeley. His research interests include data-driven decision making, with particular emphasis on addressing inefﬁciencies and inequities in health systems and physical infrastructure.

\textbf{Philip Kaminsky} Philip is the Earl J. Isaac Professor in the Science and Analysis of Decision Making in the Department of Industrial Engineering and Operations Research at UC Berkeley, where he previously served as Executive Associate Dean of the College of Engineering, faculty director of the Sutardja Center for Entrepreneurship and Technology, and department chair of Industrial Engineering and Operations Research.  His research focuses primarily on the analysis and development of robust and efficient tools and techniques for design, operation, and risk management in logistics systems and supply chains. He consults in the areas of production planning, logistics, and supply chain management., and prior to his graduate education, he worked in production engineering and control at Merck \& Co.

\textbf{Elena Flowers} Elena is an Associate Professor in the Department of Physiological Nursing and Institute for Human Genetics at the University of California, San Francisco. Her program of research is focused on the relationships between risk factors, risk-reduction interventions, and molecular biomarkers for chronic diseases, namely type-2 diabetes and preterm birth, in high risk racial groups.

\textbf{Yoshimi Fukuoka} Yoshimi, PhD, RN, FAAN is a professor in the Department of Physiological Nursing at the University of California, San Francisco. Her research focuses on the primary prevention of cardiovascular disease and type 2 diabetes using digital technologies and artificial intelligence.



\begin{APPENDICES}
\section{Proofs of Propositions in Text}
\label{ap:proof}
\proof{Proof of Proposition \ref{prop:reform_prop}: }
To obtain this formulation first we can augment the objective function of the log-likelihood problem by adding the constant term $\log \sum_{t\in\mathcal{T}_a}p(r_{a,t}|\theta^*_a,x^*_{a,t})$ and multiplying by the positive constant $\frac{1}{n(\mathcal{T}_a)}$ which does not change the value of the optimal solution. Next we use functional compositions to contract the dynamics and obtain an objective function which is explicitly a function of $\theta_a,x_{0,a}$. \halmos
\endproof

\proof{Proof of Lemma \ref{lemma:lip_lem}: } We can see that this is the case by noting that by Assumption \ref{ass:sub_gauss} we have that each of the log-likelihood ratios are Lipschitz with constant $L_p$. Since Lipschitz continuity is preserved by addition and averaging we note that the average of all of these log-likelihood ratios is also $L_p$-Lipschitz. Next we use the property that functional compositions of Lipschitz functions are Lipschitz with a constant equal to the product of their respective constants and the Lipschitz continuity is preserved through point wise maxima \citep{rockafellar2009variational}. Since the absolute value function is 1-Lipschitz and we are performing maximization we have that $\phi$ is indeed $L_p$-Lipschitz with respect to the input sequence. \halmos  
\endproof

\proof{Proof of Lemma \ref{lem:lip_param}:}
To show the first result we use a similar argument to that of the proof of Lemma \ref{lemma:lip_lem} by showing that the likelihood is Lipschitz and then using the preservation of Lipschitz continuity across functional compositions. First consider $h_a^t(x)$, since by Assumption \ref{ass:h_lip} $h_a$ is locally $L_h$ Lipschitz such that $L_h \leq 1$, hence  
 we have that with respect to $x,x' \in \mathcal{X}$ $\|h_a^t(x) - h_a^t(x')\|_2 < \|x - x'\|_2$.  Hence $h_a^t(x)$ is locally 1-Lipschitz continuous with respect to $x,t$. Next, applying Assumption \ref{ass:lip_ass} shows that since the likelihood ratio is $L_f$-Lipschitz with respect to its two inputs we simply have a composition of Lipschitz functions and the result follows. 

To show the second result note that $\ell$ depends on $t$ only through the composite dynamics mapping $h_a^t$. By definition $h^t(x) \in \mathcal{X}$ which is a bounded set, we have that for any $t,t' \in \{1,...,T\}$ $\|h_a^t(x) - h_a^{t'}(x)\|_2 \leq \diam(\mathcal{X})$, thus using Assumption \ref{ass:lip_ass} we obtain the desired result.\halmos
\endproof

\proof{Proof of Lemma \ref{lem:rademach}}
To prove this result we first bound the expectation by a Rademacher average \citep{bartlett2002} and then apply Dudley's Integral bound \citep{wainwright2015}. First let us consider the explicit form of $\mathbb{E}\varphi(\{r_{a,t}\}_{t=1}^{n(\mathcal{T}_a)})$. Using an identically distributed sequence of rewards $\{r'_{a,t}\}_{t=1}^{n(\mathcal{T}_a)}$ which is independent of the observed sequence we see that
\begin{equation}
\begin{aligned}
&\mathbb{E}\sup_{\theta_a,x_{a,0} \in \Theta\times \mathcal{X}}\big|\frac{1}{n(\mathcal{T}_a)}\sum_{t\in\mathcal{T}_a}\log\frac{ p(r_{a,t}|\theta^*_a,h_a^t(x^*_{a,0},\theta^*_{a},\pi_1^t))}{p(r_{a,t}|\theta_a,h_a^t(x_{a,0},\theta_{a},\pi_1^t))} - \frac{1}{n(\mathcal{T}_a)}D_{a,\pi_1^T}(\theta^*_a,x^*_{a,0}||\theta_a,x_{a,0}) \big|\\
&=\mathbb{E}\sup_{\theta_a,x_{a,0} \in \Theta\times \mathcal{X}}\Big|\frac{1}{n(\mathcal{T}_a)}\mathbb{E}\big[\sum_{t\in\mathcal{T}_a}\log\frac{ p(r_{a,t}|\theta^*_a,h_a^t(x^*_{a,0},\theta^*_{a},\pi_1^t))}{p(r_{a,t}|\theta_a,h_a^t(x_{a,0},\theta_{a},\pi_1^t))} - \sum_{t\in\mathcal{T}_a}\log\frac{ p(r'_{a,t}|\theta^*_a,h_a^t(x^*_{a,0},\theta^*_{a},\pi_1^t))}{p(r'_{a,t}|\theta_a,h_a^t(x_{a,0},\theta_{a},\pi_1^t))} \big|\{r_{a,t}\}_{t=1}^{n(\mathcal{T}_a)} \big] \Big|\\
&\leq \mathbb{E}\sup_{\theta_a,x_{a,0} \in \Theta\times \mathcal{X}}\Big|\frac{1}{n(\mathcal{T}_a)}\Big(\sum_{t\in\mathcal{T}_a}\log\frac{ p(r_{a,t}|\theta^*_a,h_a^t(x^*_{a,0},\theta^*_{a},\pi_1^t))}{p(r_{a,t}|\theta_a,h_a^t(x_{a,0},\theta_{a},\pi_1^t))} - \sum_{t\in\mathcal{T}_a}\log\frac{ p(r'_{a,t}|\theta^*_a,h_a^t(x^*_{a,0},\theta^*_{a},\pi_1^t))}{p(r'_{a,t}|\theta_a,h_a^t(x_{a,0},\theta_{a},\pi_1^t))}\Big) \Big|.
\end{aligned}
\end{equation}
Here the inequality follows from Jensen's Inequality \citep{qu2011}. Let $\{\epsilon_t\}_{t=1}^{n(\mathcal{T}_a)}$ be a sequence of i.i.d. Rademacher random variables, which are independent of the observations $r_{a,t},r'_{a,t}$, then through a symmetrization argument its clear that
\begin{equation} \label{eq:radmemach_avg}
\mathbb{E}\varphi(\{r_{a,t}\}_{t=1}^{n(\mathcal{T}_a)}) \leq 2\mathbb{E}\sup_{\theta_a,x_{a,0}\in \Theta \times \mathcal{X}}\Big|\frac{1}{n(\mathcal{T}_a)}\sum_{t\in\mathcal{T}_a} \epsilon_t\log\frac{ p(r_{a,t}|\theta^*_a,h_a^t(x^*_{a,0},\theta^*_{a},\pi_1^t))}{p(r_{a,t}|\theta_a,h_a^t(x_{a,0},\theta_{a},\pi_1^t))}\Big|.
\end{equation}
Since $x^*_{a,0},\theta^*_{a}$ are constants we can use simplify the above expression using the notation introduced in Lemma \ref{lem:lip_param} to $2\mathbb{E}\sup_{\theta_a,x_a \in \Theta\times\mathcal{X}}\Big|\frac{1}{n(\mathcal{T}_a)}\sum_{t\in\mathcal{T}_a} \epsilon_t\ell(\theta_a,x_{0,a},t)\Big|$.  We can bound this expression as follows
\begin{equation} \label{eq:bnd_radem}
\begin{aligned}
&2\mathbb{E}\sup_{\theta_a,x_a \in \Theta\times\mathcal{X}}\Big|\frac{1}{n(\mathcal{T}_a)}\sum_{t\in\mathcal{T}_a} \epsilon_t\ell(\theta_a,x_{0,a},t)\Big|, \\
&= 2\mathbb{E}\sup_{\theta_a,x_a \in \Theta\times\mathcal{X}}\Big|\frac{1}{n(\mathcal{T}_a)}\sum_{t\in\mathcal{T}_a} \epsilon_t(\ell(\theta_a,x_{0,a},t)-\ell(\theta_a,x_{0,a},0)+\ell(\theta_a,x_{0,a},0))\Big|, \\
&\leq 2\mathbb{E}\sup_{\theta_a,x_a \in \Theta\times\mathcal{X}}\Big|\frac{1}{n(\mathcal{T}_a)}\sum_{t\in\mathcal{T}_a} \epsilon_t(\ell(\theta_a,x_{0,a},t)-\ell(\theta_a,x_{0,a},0))\Big|+ 2\mathbb{E}\sup_{\theta_a,x_a \in \Theta\times\mathcal{X}}\Big|\frac{1}{n(\mathcal{T}_a)}\sum_{t\in\mathcal{T}_a}\epsilon_t\ell(\theta_a,x_{0,a},0)\Big|. 
\end{aligned}
\end{equation}

For our analysis we can consider each of these terms separately and bound them using Dudley's Integral Bound \citep{wainwright2015} and Lemmas \ref{lem:lip_param},\ref{lem:ent_int}. Consider the first term, note that by Lemma \ref{lem:lip_param} we have that $|\ell(\theta_a,x_{0,a},t)-\ell(\theta_a,x_{0,a},0)| \leq L_f\diam(\mathcal{X})$ and is contained in an $\ell_2$ ball of this radius, hence by Lemma \ref{lem:ent_int}
\begin{multline}
2\mathbb{E}\sup_{\theta_a,x_a \in \Theta\times\mathcal{X}}\Big|\frac{1}{n(\mathcal{T}_a)}\sum_{t\in\mathcal{T}_a} \epsilon_t(\ell(\theta_a,x_{0,a},t)-\ell(\theta_a,x_{0,a},0))\Big|, \\
\leq 8 \int_{0}^{L_f\diam(\mathcal{X})} \sqrt{\frac{\log\mathcal{N}(L_f\diam(\mathcal{X})B_2,\alpha,\|\|_2)}{n(\mathcal{T}_a)}}d\alpha \leq 8L_f\diam(\mathcal{X})\sqrt{\frac{\pi}{n(\mathcal{T}_a)}}.
\end{multline}
The last inequality follows from using  a volume bound on the covering number and using integration by parts. Next consider the second term in \eqref{eq:bnd_radem}, we can bound this term using a direct application of Dudley's entropy integral as follows
\begin{multline}\label{eq:first_integ_bnd}
2\mathbb{E}\sup_{\theta_a,x_a \in \Theta\times\mathcal{X}}\Big|\frac{1}{n(\mathcal{T}_a)}\sum_{t\in\mathcal{T}_a} \epsilon_t\ell(\theta_a,x_{0,a},0)\Big| \leq 16\sqrt{2} \int_0^{\infty} \sqrt{\frac{\log2\mathcal{N}(\alpha,\ell(\Theta\times\mathcal{X}),\|\|_2)}{n(\mathcal{T}_a)}}d\alpha, \\ \leq
 16\sqrt{2} \int_0^{\infty} \sqrt{\frac{\log2\mathcal{N}(\frac{\alpha}{L_f},\Theta\times\mathcal{X},\|\|_2)}{n(\mathcal{T}_a)}}d\alpha.
\end{multline}
Let $v_\ell B_2$ be the $\ell_2$ ball on $\mathbb{R}^{d_x+d_\theta}$ with radius $v_\ell = \diam(\mathcal{X}\times\Theta)$, then
\begin{equation} \label{eq:integ_bound}
\eqref{eq:first_integ_bnd} \leq 16\sqrt{2} \int_0^{\infty} \sqrt{\frac{\log2\mathcal{N}(\frac{\alpha}{L_f},B_\ell,\|\|_2)}{n(\mathcal{T}_a)}}d\alpha \leq 16\sqrt{2} \int_0^{\infty} \sqrt{\frac{\log2(\frac{3v_\ell L_f}{\alpha})^{d_x+d_\theta}}{n(\mathcal{T}_a)}}d\alpha 
\end{equation}
Solving the integral shows that \eqref{eq:integ_bound} $\leq 48\sqrt{2}(2)^{\frac{1}{d_x+d_\theta}}L_fv_\ell\sqrt{\frac{\pi(d_x+d_\theta)}{n(\mathcal{T}_a)}} $. Hence the result follows. \halmos
\endproof

\proof{Proof of Theorem \ref{thm:two_side_ineq}:} 
Lemma \ref{lemma:lip_lem} guarantees that the mapping $\varphi$ is Lispschitz continuous with respect to the observed rewards with parameter $L_p$, furthermore we have by Assumption \ref{ass:sub_gauss} that the reward distributions are sub-Gaussian with parameter $\sigma^2$. By applying Theorem 1 from \cite{kontorovich2014} we obtain for $\xi >0$:
\begin{equation}
\mathbb{P}\Big(\varphi(\{r_t\}_{t=1}^{n(\mathcal{T}_a)}) - \mathbb{E}\varphi(\{r_t\}_{t=1}^{n(\mathcal{T}_a)}) > \xi \Big) \leq \exp(\frac{-\xi^2n(\mathcal{T}_a)}{2L_p^2\sigma^2}).
\end{equation}

Hence, using the upper bound obtained from Lemma \ref{lem:rademach}, we can substitute the result into the above equation giving the desired result. \halmos
\endproof

\proof{Proof of Theorem \ref{thm:concent_eq}: } Using Theorem \ref{thm:two_side_ineq} we know that with probability at least $1-\exp(\frac{-\xi^2n(\mathcal{T}_a)}{2L^2_p\sigma^2})$ we have:
\begin{equation}
\frac{1}{n(\mathcal{T}_a)}D_{a,\pi_1^T}(\theta^*_a,x^*_{a,0}||\hat{\theta}_a,\hat{x}_{a,0}) - \frac{1}{n(\mathcal{T}_a)}\sum_{t\in n(\mathcal{T}_a)} \log \frac{p(r_{a,t}|\theta_a^*,h_a^t(x^*_{a,0},\theta^*_a,\pi_1^t))}{p(r_{a,t}|\hat{\theta}_a,h_a^t(\hat{x}_{a,0},\hat{\theta}_a,\pi_1^t))} \leq \frac{c_f(d_x,d_\theta)}{\sqrt{n(\mathcal{T}_a)}} +\xi.
\end{equation}
Also since $\hat{\theta}_a,\hat{x}_a$ are minimizers of the empirical trajectory divergence implies that
\begin{equation}
\frac{1}{n(\mathcal{T}_a)}\sum_{t\in n(\mathcal{T}_a)} \log \frac{p(r_{a,t}|\theta_a^*,h_a^t(x^*_{a,0},\theta^*_a,\pi_1^t))}{p(r_{a,t}|\hat{\theta}_a,h_a^t(\hat{x}_{a,0},\hat{\theta}_a,\pi_1^t))} \leq \frac{1}{n(\mathcal{T}_a)}\sum_{t\in n(\mathcal{T}_a)} \log \frac{p(r_{a,t}|\theta_a^*,h_a^t(x^*_{a,0},\theta^*_a,\pi_1^t))}{p(r_{a,t}|\theta^*_a,h_a^t(x^*_{a,0},\theta^*_a,\pi_1^t))} = 0.
\end{equation}
Hence the desired result follows. \halmos
\endproof

\proof{Proof of Proposition \ref{prop:init_reg_bnd}: } Recall that by definition $\mathbb{E}R_\Pi(T) = \sum_{t=1}^Tg(\theta_{pi^*_t},x_{pi^*_t}) - g(\theta_{pi_t},x_{pi_t})$. Since by Assumption \ref{ass:lip_ass} we have that $g$ is $L_g$-Lipschitz then we have $\forall t$ that $g(\theta_{pi^*_t},x_{pi^*_t}) - g(\theta_{pi_t},x_{pi_t}) \leq L_g\|(\theta_{pi^*_t},x_{pi^*_t}) - (\theta_{pi_t},x_{pi_t})\| \leq L_g\diam(\mathcal{X}\times \Theta)\mathbb{P}(\pi_t \neq \pi^*t)$. Hence  
\begin{equation}
\begin{aligned}
\mathbb{E}R_\Pi(T) &\leq L_g \diam(\mathcal{X}\times\Theta) \sum_{t=0}^T \mathbb{P}(\pi_t \neq \pi_t^*) = L_g \diam(\mathcal{X}\times\Theta) \sum_{t=0}^T \sum_{a\in\mathcal{A}}\mathbb{P}(\pi_t = a, a \neq \pi_t^*) \\
&= L_g \diam(\mathcal{X}\times\Theta)  \sum_{a\in\mathcal{A}} \sum_{t=0}^T \mathbb{P}(\pi_t = a, a \neq \pi_t^*) = L_g \diam(\mathcal{X}\times\Theta) \sum_{a\in\mathcal{A}} \mathbb{E}\tilde{T}_a  \halmos
\end{aligned}
\end{equation}
\endproof

\proof{Proof of Proposition \ref{prop:wrong_pulls}: } We proceed to prove this proposition in a similar method to that presented in \cite{auer2002}. Suppose that at time $t$, the ROGUE-UCB policy chooses $a \neq \pi^*_t$. If the upper confidence bounds hold then we observe that $g_{a,t}^{UCB} \geq g_{\pi^*_t,t}^{UCB} \geq g_{\pi^*_t,t}$. Also define the mapping $\psi_a(\gamma) = \max\{|g(\theta,h_a^t(x_{0})) - g(\hat{\theta}_a,h_a^t(\hat{x}_{a,0}))|: \frac{1}{n(\mathcal{T}_a)} D_{a,\pi_1^T}(\theta,x_0||\hat{\theta}_a,\hat{x}_{a,0}) \leq \gamma \}$. Then clearly $g_{a,t}^UCB - g(\hat{\theta}_a,h_a^t(\hat{x}_{a,0})) \leq \psi_a(A(t)\sqrt{\frac{4\log(t)}{n(\mathcal{T}_a)}})$ and $g(\hat{\theta}_a,h_a^t(\hat{x}_{a,0})) - g_{a,t} \leq \psi_a(A(t)\sqrt{\frac{4\log(t)}{n(\mathcal{T}_a)}})$. Hence we have that $g_{a,t}^{UCB} \leq 2\psi(A(t)\sqrt{\frac{4\log(t)}{n(\mathcal{T}_a)}}) + g_{a,t}$. Therefore $\psi(A(t)\sqrt{\frac{4\log(t)}{n(\mathcal{T}_a)}}) \geq \frac{1}{2}(g_{\pi^*_t,t} - g_{a,t})$. By definition of $\epsilon_a$ we thus have that $\psi(A(t)\sqrt{\frac{4\log(t)}{n(\mathcal{T}_a)}}) \geq \frac{\epsilon_a}{2}$. Therefore, by definition of $\delta_a$ we observe that $A(t)\sqrt{\frac{4\log t}{n(\mathcal{T}_a)}} \geq \delta_a$ and hence $n(\mathcal{T}_a) \leq \frac{4A(t)^2\log t}{\delta_a^2}$.

Now, consider $\tilde{T}_a$:
\begin{align}
&\tilde{T}_a = \sum_{t=1}^T \mathbf{1}\{\pi_t = a, a \neq \pi^*_t\}\\
 & = \sum_{t=1}^{T}\mathbf{1}\{\pi_t = a, a \neq \pi^*_t, n(\mathcal{T}_a) \leq \frac{4A(t)^2\log t}{\delta_a^2}\} + \sum_{t=1}^{T}\mathbf{1}\{\pi_t = a, a \neq \pi^*_t, n(\mathcal{T}_a) > \frac{4A(t)^2\log t}{\delta_a^2}\} \\
&\leq \sum_{t=1}^{T} \mathbf{1}\{\pi_t = a, a \neq \pi^*_t, n(\mathcal{T}_a) \leq \frac{4A(|\mathcal{A}|)^2\log T}{\delta_a^2}\} +\sum_{t=1}^{T}\mathbf{1}\{\pi_t = a, a \neq \pi^*_t, n(\mathcal{T}_a) > \frac{4A(t)^2\log t}{\delta_a^2}\} \\
& \leq \frac{4\log(T)}{\delta_a^2} A(|\mathcal{A}|)^2 + \sum_{t=1}^{T}\mathbf{1}\{\pi_t = a, a \neq \pi^*_t, n(\mathcal{T}_a) > \frac{4A(t)^2\log t}{\delta_a^2}\}
\end{align}
Observe that if we play sub optimal action $a$ at time $t$ this means we either severely over estimate the value of $g_{a,t}$, severely under estimate the value of $g_{\pi_t^*,t}$, or the two values are very close to each other. Hence
 \begin{multline}
 \{\pi_t=a, a\neq\pi^*_t, n(\mathcal{T}_a)  > \frac{4A(t)^2\log t}{\delta_a^2}\} \subseteq \underbrace{\{g_{a,t}^{UCB} - g_{a,t} > 2 \psi_a(A(t)\sqrt{\frac{4\log(t)}{n(\mathcal{T}_a)}}), n(\mathcal{T}_a) >\frac{4A(t)^2\log t}{\delta_a^2} \}}_{(a)} \\ \cup \underbrace{\{g_{\pi^*_t,t} >g^{\pi_t^*,t}_{UCB}, n(\mathcal{T}_a) >\frac{4A(t)^2\log t}{\delta_a^2}\}}_{(b)}\cup\underbrace{\{g_{\pi^*_t,t} - g_{a,t} \leq 2\psi_a(A(t)\sqrt{\frac{4\log(t)}{n(\mathcal{T}_a)}}), n(\mathcal{T}_a) >\frac{4A(t)^2\log t}{\delta_a^2}\}}_{(c)}.
\end{multline} 
 However, as we established in the beginning of the proof the event $(c) = \emptyset$. Also note that for events $(a),(b)$ to occur this would imply that $\theta_a,x_{a,0}$ and $\theta_{\pi_t^*},x_{\pi_t^*,0}$ are not feasible points of their respective UCB deriving problems, hence
 \begin{equation}
 \begin{aligned}
 &\{\pi_t=a, a\neq\pi^*_t, n(\mathcal{T}_a)  > \frac{4A(t)^2\log t}{\delta_a^2}\} \subseteq \{\exists s<t :\;  \frac{1}{s}D_{\pi^*_t,\pi_1^s}(\hat{\theta}_{\pi^*_t},\hat{x}_{\pi^*_t,0}||\theta_{\pi^*_t},x_{\pi^*_t,0}) > A(t)\sqrt{\frac{4\log(t)}{s}}\} \\
  & \cup\{\exists s'<t: \;\frac{1}{s'}D_{a,\pi_1^{s'}}(\hat{\theta}_{a},\hat{x}_{a,0}||\theta_{a},x_{a,0}) > A(t)\sqrt{\frac{4\log(t)}{s'}}\} \\
 &\subseteq
 \bigcup_{1\leq s < t} \{ \frac{1}{s}D_{\pi^*_t,\pi_1^s}(\hat{\theta}_{\pi^*_t},\hat{x}_{\pi^*_t,0}||\theta_{\pi^*_t},x_{\pi^*_t,0}) > A(t)\sqrt{\frac{4\log(t)}{s}}\} \bigcup_{1\leq s' < t} \{\frac{1}{s'}D_{a,\pi_1^{s'}}(\hat{\theta}_{a},\hat{x}_{a,0}||\theta_{a},x_{a,0}) > A(t)\sqrt{\frac{4\log(t)}{s'}}\}.  
 \end{aligned}
 \end{equation} 
 Taking the expected value of $\tilde{T}_a$ we obtain
 \begin{equation}
 \begin{aligned}
 &\mathbb{E}\tilde{T}_a \leq \frac{4\log(T)}{\delta_a^2} A(|\mathcal{A}|)^2 + \mathbb{E}\sum_{t=1}^{T}\mathbf{1}\{\pi_t = a, a \neq \pi^*_t, n(\mathcal{T}_a) > \frac{4A(t)^2\log t}{\delta_a^2}\}
 \\
 & \leq \frac{4\log(T)}{\delta_a^2} A(|\mathcal{A}|)^2 + \sum_{t=1}^{T}\sum_{s=1}^t\sum_{s'=1}^t\mathbb{P}(\frac{1}{s}D_{\pi^*_t,\pi_1^s}(\hat{\theta}_{\pi^*_t},\hat{x}_{\pi^*_t,0}||\theta_{\pi^*_t},x_{\pi^*_t,0}) > A(t)\sqrt{\frac{4\log(t)}{s}})\\ 
 &+ \sum_{t=1}^{T}\sum_{s=1}^t\sum_{s'=1}^t\mathbb{P}(\frac{1}{s'}D_{a,\pi_1^{s'}}(\hat{\theta}_{a},\hat{x}_{a,0}||\theta_{a},x_{a,0}) > A(t)\sqrt{\frac{4\log(t)}{s'}})  \\
 & \leq \frac{4\log(T)}{\delta_a^2} A(|\mathcal{A}|)^2 +  2\sum_{t=1}^T\sum_{s=1}^t\sum_{s'=1}^t t^{-4} \leq \frac{4\log(T)}{\delta_a^2} A(|\mathcal{A}|)^2 + \frac{\pi^2}{3}.
 \end{aligned}
 \end{equation}
 Here the third inequality is derived by Theorem \ref{thm:concent_eq} and the final inequality by utilizing the solution to the Basel Problem \citep{rockafellar2009variational}. Hence we obtain the desired result. \halmos
\endproof

\proof{Proof of Theorem \ref{thm:regret_them}: }
Using Proposition \ref{prop:init_reg_bnd} we bound the expected regret as $\mathbb{E}R_\Pi(T) \leq L_g\diam(\mathcal{X}\times\Theta)\sum_{a\in\mathcal{A}} \mathbb{E}\tilde{T}_a$. Then applying the result of Proposition \ref{prop:wrong_pulls} we obtain the desired result. \halmos
\endproof

\proof{Proof of Corollary \ref{cor:parameter_free}: }
To obtain the desired rate we must first express the bound from Theorem \ref{thm:regret_them} in terms of $\epsilon$. We can do this by noting that from Assumption \ref{ass:sub_gauss} we know that the reward distributions are sub-Gaussian with parameter $\sigma$, therefore using the property that exponential families achieve the maximum entropy given moment constraints \citep{bickel2006}, we observe that $\forall \theta'_a,x'_{a,0} \in \Theta \times \mathcal{X}$: 

	\begin{equation}
	\frac{1}{n(\mathcal{T}_a)}D_{a,\pi_1^T}(\theta_a,x_{a,0}||\theta'_a,x'_{a,0}) \geq \frac{1}{2\sigma n(\mathcal{T}_a)} \sum_{t \in \mathcal{T}}(g(\theta_a,x_{a,t})- g(\theta'_a,x'_{a,t}))^2 \geq \frac{1}{2\sigma}\epsilon_a^2 
	\end{equation}

Thus applying the definition of $\delta_a$ we observe that in fact $\delta_a \geq\frac{1}{2\sigma}\epsilon_a^2 $. Hence substituting into the bound from Theorem \ref{thm:regret_them} we obtain the following weaker regret bound in terms of the $\epsilon_a$: 

	\begin{equation}
	\mathbb{E}R_\Pi(T) \leq L_g\diam(\mathcal{X}\times\Theta)\sum_{a\in\mathcal{A}} \Bigg(A(|\mathcal{A}|)^2 \frac{16\sigma^2\log T}{\epsilon_a^4} + \frac{\pi^2}{3}\Bigg).
	\end{equation}

Note that this bound is convex and monotonically decreasing in the $\epsilon_a$, and in fact is minimized if $\epsilon_a = \max_{(\theta_1,x_1),(\theta_2,x_2) \in \Theta\times\mathcal{X}} |g(\theta_1,x_1) - g(\theta_2,x_2)|$ for all $a \in \mathcal{A}$. Hence all $\epsilon_a$ will be equal to each other when this bound is minimized. However, this bound will go to infinity as the $\epsilon_a$ become small. Note though that since by assumption, $\epsilon_a \geq \epsilon$ and hence as $\epsilon_a$ become small we would expect our regret to be at a rate of $\epsilon T$. Using these two notions we can thus conclude that our regret can be bounded as:

	\begin{equation} \label{eq:min_bnd}
	\mathbb{E}R_\Pi(T) \leq \min\{ L_g\diam(\mathcal{X}\times\Theta)|\mathcal{A}| \Big(A(|\mathcal{A}|)^2 \frac{16\sigma^2\log T}{\epsilon^4} + \frac{\pi^2}{3}\Big),\epsilon T\}
	\end{equation}

Since one of these expressions increases in $\epsilon$ and the other decreases in $\epsilon$ we can minimize the right hand side by setting both parts of the $\min$ to be equal to each other. This shows that the minimizing epsilon can be given as a zero to the following quintic polynomial:

	\begin{equation}
	T\epsilon^5 - L_g\diam(\mathcal{X}\times\Theta)|\mathcal{A}|\frac{\pi^2}{3}\epsilon^4 - L_g\diam(\mathcal{X}\times\Theta)|\mathcal{A}|A(|\mathcal{A}|)^216\sigma^2\log T
	\end{equation}

Since the zeros of this polynomial are difficult to compute we can use the Fujiwara bound on polynomial zeros \citep{fujiwara1916} to show that for sufficiently large $T$, $\epsilon = \mathcal{O}\Bigg(\Bigg( \frac{|\mathcal{A}|A(|\mathcal{A}|)^2\sigma^2\log T}{T}\Bigg)^{\frac{1}{5}}\Bigg)$, and hence the desired result follows from substituting this value into Equation \eqref{eq:min_bnd}. \halmos

\endproof

\proof{Proof of Theorem \ref{thm:eps_greed_reg}: }
We proceed to prove this theorem using similar analysis to that provided in \cite{auer2002finite}. The probability that an arm $a$ is played sub-optimally by the $\epsilon$-ROGUE algorithm is given by:
\begin{equation}
\mathbb{P}(\pi_t=a,a\neq\pi^*_t) = \frac{\epsilon_t}{|\mathcal{A}|} + (1-\epsilon_t)\mathbb{P}\big(g(\hat{x}_{a,t},\hat{\theta}_a) \geq g(\hat{x}_{\pi_t^*,t},\hat{\theta}_{\pi_t^*})\big)
\end{equation}
Where the first term coms from the probability of performing an exploration step and the second term corresponds to an exploration step. First lets consider the exploitation term, define $\Delta_{a,t} = g(x_{\pi_t^*,t},\theta_{\pi^*_t}) - g(x_{a,t},\theta_a)$, using this value we can upper bound the exploitation term as follows:

\begin{equation}
\begin{aligned}
&\mathbb{P}\big(g(\hat{x}_{a,t},\hat{\theta}_a) \geq g(\hat{x}^*_{a,t},\hat{\theta}^*_a)\big)  \\
  &\leq \mathbb{P}\big( (g(\hat{x}_{a,t},\hat{\theta}_{a}) - g(x_{a,t},\theta_{a})) - (g(\hat{x}_{\pi*_t,t},\hat{\theta}_{\pi^*_t}) - g(x_{\pi^*_t,t},\theta_{\pi^*_t})) \geq \Delta_{a,t} \big) \\
  &\leq \mathbb{P}\big(g(\hat{x}_{a,t},\hat{\theta}_{a}) - g(x_{a,t},\theta_{a}) \geq \frac{\Delta_{a,t}}{2} \big) + \mathbb{P}\big( g(\hat{x}_{\pi*_t,t},\hat{\theta}_{\pi^*_t}) - g(x_{\pi^*_t,t},\theta_{\pi^*_t}) \leq \frac{-\Delta_{a,t}}{2} \big) \\
\end{aligned}
\end{equation}

Using the definitions of $\epsilon_a,\delta_a$ from Theorem \ref{thm:regret_them} we obtain the following:
\begin{equation}
\begin{aligned}
&\leq \mathbb{P}\big(g(\hat{x}_{a,t},\hat{\theta}_{a}) - g(x_{a,t},\theta_{a}) \geq \frac{\epsilon_a}{2} \big) + \mathbb{P}\big( g(\hat{x}_{\pi*_t,t},\hat{\theta}_{\pi^*_t}) - g(x_{\pi^*_t,t},\theta_{\pi^*_t}) \leq \frac{-\epsilon_a}{2} \big) \\
&\leq \mathbb{P}\big(\bar{D}_{a,\pi_1^t}(\theta_a,x_{a,0}||\hat{\theta}_a,\hat{x}_{a,0}) \geq \delta_a \big) + \mathbb{P}\big( \bar{D}_{\pi^*_t,\pi_1^t}(\theta_{\pi^*_t},x_{\pi^*_t,0}||\hat{\theta}_{\pi^*_t},\hat{x}_{\pi^*_t,0}) \geq \delta_a \big)
\end{aligned}
\end{equation}

Both of these terms can be bounded in a similar manner, as such we will present the arguments for bounding the first term and note that a similar procedure can be used to bound the second. We begin as follows:
\begin{equation}
\begin{aligned}
&\mathbb{P}\big(\bar{D}_{a,\pi_1^t}(\theta_a,x_{a,0}||\hat{\theta}_a,\hat{x}_{a,0}) \geq \delta_a \big) \leq \sum_{n=1}^t \mathbb{P}\big(n(\mathcal{T}_a)=n, \bar{D}_{a,\pi_1^n}(\theta_a,x_{a,0}||\hat{\theta}_a,\hat{x}_{a,0}) \geq \delta_a \big) \\
& \sum_{n=1}^t \mathbb{P}\big(n(\mathcal{T}_a)=n| \bar{D}_{a,\pi_1^n}(\theta_a,x_{a,0}||\hat{\theta}_a,\hat{x}_{a,0}) \geq \delta_a \big)\mathbb{P}(\bar{D}_{a,\pi_1^n}(\theta_a,x_{a,0}||\hat{\theta}_a,\hat{x}_{a,0}) \geq \delta_a)
\end{aligned}
\end{equation}

Applying Theorem \ref{thm:concent_eq} we obtain:
\begin{equation}
\begin{aligned}
& \sum_{n=1}^t \mathbb{P}\big(n(\mathcal{T}_a)=n| \bar{D}_{a,\pi_1^n}(\theta_a,x_{a,0}||\hat{\theta}_a,\hat{x}_{a,0}) \geq \delta_a \big)\exp\big(\frac{-(\delta_a\sqrt{n}-c_f(d_x,d_\theta))^2}{2L_p^2\sigma^2}\big) \\
& \leq \sum_{n=1}^{\lfloor w \rfloor}\mathbb{P}\big(n(\mathcal{T}_a)=n| \bar{D}_{a,\pi_1^n}(\theta_a,x_{a,0}||\hat{\theta}_a,\hat{x}_{a,0}) \geq \delta_a \big) + \sum_{n= \lceil w \rceil}^\infty \exp\big(\frac{-(\delta_a\sqrt{n}-c_f(d_x,d_\theta))^2}{2L_p^2\sigma^2}\big) \\
\end{aligned}
\end{equation}

Where we define $w = \frac{1}{2|\mathcal{A}|}\sum_{n=1}^t \epsilon_t$. We can bound the first term in the expression above by using the the total number of times arm $a$ is chosen as part of an exploration step, we call this quantity $n^R(\mathcal{T}_a)$. This can be done as follows:
\begin{equation}
\begin{aligned}
&\sum_{n=1}^{\lfloor w \rfloor}\mathbb{P}\big(n(\mathcal{T}_a)=n| \bar{D}_{a,\pi_1^n}(\theta_a,x_{a,0}||\hat{\theta}_a,\hat{x}_{a,0}) \geq \delta_a \big) \\
&\leq \sum_{n=1}^{\lfloor w \rfloor} \mathbb{P}\big(n^R(\mathcal{T}_a)\leq n |\bar{D}_{a,\pi_1^n}(\theta_a,x_{a,0}||\hat{\theta}_a,\hat{x}_{a,0}) \geq \delta_a\big) \\
&\leq w\mathbb{P}\big(n^R(\mathcal{T}_a)\leq w\big)
\end{aligned}
\end{equation}

The last two inequality follows since the amount of exploration steps is independent of how the parameter estimates are concentrating. Note that $w = \frac{1}{2}\mathbb{E}n^R(\mathcal{T}_a)$, and that $\Var(n^R(\mathcal{T}_a)) \leq \mathbb{E}n^R(\mathcal{T}_a)$. Hence we can employ the Bernstein bound \citep{wainwright2015} to obtain the following upper bound on the above expression:
\begin{equation}
w\mathbb{P}\big(n^R(\mathcal{T}_a)\leq w\big) \leq w \exp\big(\frac{-3w}{14}\big)
\end{equation}

 We will bound the second term in the summation using an integral bound. To do this, first note that since $-(\delta_a\sqrt{n}-c_f(d_x,d_\theta))^2$ is concave in $n$, so we can upper bound it using an appropriate first order approximation. Namely, consider the first order Taylor expansion about $n = (\frac{2c_f(d_\theta,d_x)}{\delta_a})^2$, then:
\begin{equation}
	\begin{aligned}
	& -(\delta_a\sqrt{n} - c_f(d_x,d_\theta))^2 \leq \frac{-\delta_a^2}{2}n + c_f(d_x,d_\theta)^2 \\
	&\implies \exp\big(\frac{-(\delta_a\sqrt{n} - c_f(d_x,d_\theta))^2}{2L_p^2\sigma^2}\big) \leq \exp\big(\frac{\frac{-\delta_a^2}{2}n + c_f(d_x,d_\theta)^2}{2L_p^2\sigma^2}\big)
	\end{aligned}
\end{equation} 
Then using this bound we can upper bound the second summation term as follows:
\begin{equation}
\begin{aligned}
&\sum_{n= \lceil w \rceil}^\infty \exp\big(\frac{-(\delta_a\sqrt{n}-c_f(d_x,d_\theta))^2}{2L_p^2\sigma^2}\big) \leq \sum_{n= \lceil w \rceil}^\infty \exp\big(\frac{\frac{-\delta_a^2}{2}n + c_f(d_x,d_\theta)^2}{2L_p^2\sigma^2}\big) \\
& \leq \frac{4L_p^2\sigma^2}{\delta_a^2}\exp\big(\frac{\frac{-\delta_a^2}{2}\lfloor w \rfloor + c_f(d_x,d_\theta)^2}{2L_p^2\sigma^2}\big)
\end{aligned}
\end{equation}

Combining all of these equations together we obtain that:
\begin{equation}\label{eq:prob_w_bound}
\mathbb{P}(\pi_t=a,a\neq\pi^*_t) \leq \frac{\epsilon_t}{|\mathcal{A}|} + 2w \exp\big(\frac{-3w}{14}\big) + \frac{8L_p^2\sigma^2}{\delta_a^2}\exp\big(\frac{\frac{-\delta_a^2}{2}\lfloor w \rfloor + c_f(d_x,d_\theta)^2}{2L_p^2\sigma^2}\big)
\end{equation}

All that remains to complete the proof is to provide a lower bound on $w$. Let $t': = \frac{c|\mathcal{A}|}{\delta_{min}^2}$ then:
\begin{equation}
\begin{aligned}
&w = \frac{1}{2|\mathcal{A}|}\sum_{n=1}^t \epsilon_n = \frac{1}{2|\mathcal{A}|}\sum_{n=1}^{t'} \epsilon_n + \frac{1}{2|\mathcal{A}|}\sum_{n=t'}^t \epsilon_n \\
& = \frac{c}{2\delta_{min}^2} + \frac{c}{2\delta_{min}^2}\sum_{n=t'}^{t} \frac{1}{n} \\
& \geq \frac{c}{2\delta_{min}^2} + \frac{c}{2\delta_{min}^2}\log \frac{t\delta_{min}^2}{c|\mathcal{A}|} = \frac{c}{2\delta_{min}^2}\log \frac{et\delta_{min}^2}{c|\mathcal{A}|}
\end{aligned}
\end{equation}

Substituting thins bound into \eqref{eq:prob_w_bound} and expanding $\epsilon_t$ gives the desired result. \halmos

\endproof

\section{Technical Metric Entropy Lemma}
\begin{lemma} \label{lem:ent_int}
Let $a \in A \subseteq \mathbb{R}^n$ such that $A$ is bounded and $K = \max_{a\in A} \frac{d(a,0)}{n}$ with respect to some metric $d$ and $\forall a \in A, \|a\|_2 \leq d(a,0)$ . Then for i.i.d Rademacher process $\{\epsilon_i\}_{i=1}^n$ :
\begin{equation}
\mathbb{E} \sup_{a\in A} | \frac{1}{n} \sum_{i=1}^{n} \epsilon_i a_i | \leq 4 \int_{0}^{K} \sqrt{\frac{\log2\mathcal{N}(\alpha,A,d)}{n}} d\alpha
\end{equation}
\end{lemma}
\proof{Proof: }We proceed to prove this result in a similar technique to that used by \cite{wainwright2015}. Let  $\bar{A} = A \cup A^-$ and $\{\hat{A}_i\}_{i=0}^N$ be a sequence of successively finer covers of set $\bar{A}$, such that $\hat{A}_i$ is an $\alpha_i$ cover of set $\bar{A}$ with respect to metric $d$ and $\alpha_i = 2^{-i}K$. Next, define a sequence of  approximating vectors of $a$ and denote these by $\hat{a}_i$ such that for any two successive approximations $\hat{a}_i \in \hat{A}_i$ and $\hat{a}_{i-1} \in \hat{A}_{i-1}$ we have that $d(\hat{a}_i,\hat{a}_{i-1}) \leq \alpha_i$. Then observe we can rewrite $a$ as follows:
\begin{equation}
a = a + \hat{a}_N - \hat{a}_N = \hat{a}_0 + \sum_{i=1}^{N}(\hat{a}_i - \hat{a}_{i-1}) + a - \hat{a}_n
\end{equation}
Observe that we can set $\hat{a}_0$ to the 0 vector since clearly a metric ball of radius $K$ will form a $K$ cover of set $A$. Hence we obtain:
\begin{align}
&\mathbb{E} \sup_{a\in A} | \frac{1}{n} \sum_{i=1}^{n}\epsilon_i a_i |= \mathbb{E}\sup_{a\in A} | \frac{1}{n} \langle \epsilon, a \rangle | = \mathbb{E} \sup_{a \in \bar{A}} \frac{1}{n} \langle \epsilon, a \rangle = \mathbb{E}\sup_{a\in \bar{A}}  \frac{1}{n} \langle \epsilon, \sum_{j=1}^{N}(\hat{a}_i - \hat{a}_{i-1}) + a - \hat{a}_N \rangle \\ 
& \leq \mathbb{E} \sum_{j=1}^{N}\sup_{\hat{a}_j \in \hat{A}_j, \hat{a}_{j-1} \in \hat{A}_{j-1}} \langle \epsilon, \hat{a}_j - \hat{a}_{j-1}\rangle + \mathbb{E} \sup_{a\in \bar{A}} \langle \epsilon, a - \hat{a}_N \rangle  \\
&\leq \sum_{j=1}^N \alpha_i \sqrt{\frac{2 \log |\hat{A}_j||\hat{A}_{j-1}|}{n}} + \alpha_N \label{eq:finit_class}
\end{align}
Here the final inequality is obtained by applying the finite class lemma \citep{wainwright2015}. Observe that $|\hat{A}_{j-1}| \leq |\hat{A}_{j-1}| = \mathcal{N}(\alpha_i,\bar{A},d)$ and that by construction $\alpha_j = 2(\alpha_j - \alpha_{j+1})$. Hence:
\begin{align}
&\mathbb{E} \sup_{a\in A} | \frac{1}{n} \sum_{i=1}^{n}\epsilon_i a_i | \leq \sum_{j=1}^{N} 4(\alpha_j - \alpha_{j+1}) \sqrt{\frac{\log \mathcal{N}(\alpha_i,\bar{A},d)}{n}} + \alpha_N \\ 
&\leq 4 \int_{\alpha_{N+1}}^{\alpha_0} \sqrt{\frac{\log \mathcal{N}(\alpha,\bar{A},d)}{n}} d\alpha + \alpha_N \rightarrow 4 \int_0^K \sqrt{\frac{\log \mathcal{N}(\alpha,\bar{A},d)}{n}} d\alpha
\end{align}
Note that $\mathcal{N}(\alpha,\bar{A},d) \leq 2\mathcal{N}(\alpha,A,d)$ thus completing the proof. \halmos
\endproof

\end{APPENDICES}
\end{document}